\documentclass[12pt]{amsart}
\usepackage{amsmath}
\usepackage{amsthm}
\usepackage{amssymb}
\usepackage{amscd}
\usepackage{amsfonts}
\usepackage{amsbsy}
\usepackage{epsfig,afterpage}
\usepackage[dvips]{psfrag}
\usepackage{subfigure}
\usepackage{srcltx}

\usepackage{verbatim}
\usepackage{graphicx}

\usepackage{multirow}
\usepackage{graphics}
\usepackage{epstopdf}
\usepackage{overpic}
\usepackage{multicol}

\usepackage[colorlinks=true, linkcolor=blue, citecolor=red, urlcolor=blue]{hyperref}

\newtheorem {theorem} {Theorem}%[section]
\newtheorem {proposition} [theorem]{Proposition}

\newtheorem {lemma}  [theorem]{Lemma}

\newtheorem {remark} [theorem]{Remark}

\newtheorem {definition} [theorem]{Definition}

\newtheorem{mtheorem}{Theorem}

\usepackage{float}
\usepackage{tikz, wrapfig}
\tikzset{node distance=3cm, auto}
\allowdisplaybreaks

\begin{document}

\title[Generalized upper principal part]
{Generalized upper principal part of real planar polynomial vector fields}

\author[T. M. Dalbelo, R. Oliveira and O. H. Perez]
{Tha\'is Maria Dalbelo$^{1}$, Regilene Oliveira$^{2}$, Otavio Henrique Perez$^{2}$}

\address{$^{1}$Federal University of S\~{a}o Carlos (UFSCar). Rodovia Washington Luís, Km 235, Zip Code 13565-905, S\~{a}o Carlos, S\~{a}o Paulo, Brazil.}

\address{$^{2}$University of S\~{a}o Paulo (USP), Institute of Mathematics and Computer Science. Avenida Trabalhador S\~{a}o Carlense, 400, Zip Code 13566-590, S\~{a}o Carlos, S\~{a}o Paulo, Brazil.}

\email{thaisdalbelo@ufscar.br}
\email{regilene@icmc.usp.br}
\email{otavio.perez@icmc.usp.br}

\thanks{ .}


\date{}
\dedicatory{}

\begin{abstract}
The goal of this paper is to generalize recent results about the topological classification of the dynamics of a real planar polynomial vector field near infinity. Given a polynomial vector field $X$, using Newton polyhedra one can define its generalized upper principal part $X_{\Gamma}^{U}$. By dropping some monomials of $X_{\Gamma}^{U}$, we define its minimal generalized upper principal part $X_{G}^{U}$. We prove that there exist an open and dense set $\widetilde{\mathfrak{U}}_{1}$ in the set of polynomial vector fields with Newton degenerate upper principal part satisfying the following property: if $X\in\widetilde{\mathfrak{U}}_{1}$, then $X$ and $X_{G}^{U}$ are topologically equivalent near infinity, provided that $X_{G}^{U}$ satisfies some additional non-degeneracy assumptions. Our techniques rely on toric compactification and the Normal Form Theorem.
\end{abstract}

\maketitle

\section{Introduction}

An important problem in the study of vector fields near singular points is to determine their topological type by reducing the original system to a simpler one. In this direction, the Newton polytope is a useful tool for describing the dynamics of planar vector fields. For instance, if $0\in\mathbb{R}^{2}$ is an isolated non-monodromic singularity of an analytic vector field $X$ and its \textit{principal part} $X_{\Delta}^{L}$ is Newton non-degenerate, then $X_{\Delta}^{L}$ is topologically equivalent to $X$ near the origin (see \cite{Ber78,BM90}). This means that there exists a homeomorphism that maps phase curves of $X_{\Delta}^{L}$ to phase curves of $X$, preserving orientation (but not necessarily preserving time). We refer to \cite{Alo2015} for an analogous result in dimension three. Throughout this paper, the principal part of $X$ addressed in \cite{Ber78,BM90} will be called the \textit{lower principal part} (see Section \ref{sec-lower} for precise definitions and a statement of this result). The existence of a topological equivalence between the analytic vector field $X$ and its lower principal part $X_{\Delta}^{L}$ in the Newton degenerate case was studied in \cite{Zup96,Zup00}.

In \cite{Pel95} and \cite{Pan06}, an algorithm for the resolution of singularities of analytic vector fields in dimensions two and three, respectively, was given. In both papers, the techniques relied on the use of the Newton polygon. An algorithm for the resolution of singularities (also based on the Newton polygon) of real analytic constrained differential systems was given in \cite{PerSil22-ii}. In \cite{Bru2000}, the author studied the asymptotics of integral curves of vector fields near singular points and near infinity. Problems related to the Newton polygon and integrability of vector fields can be found in \cite{DGV22} and the references therein. Finally, we refer to \cite{AGR2011,AGR2014} for problems related to monodromic singularities and Newton polygons.

Recently, the use of the Newton polytope in the study of the topological classification of the dynamics of a polynomial vector field near infinity was addressed in \cite{DOP24, OliVal26}. More precisely, in \cite{DOP24}, the authors proved that, if the \textit{upper principal part} $X_{\Delta}^{U}$ of a real planar polynomial vector field $X$ is Newton non-degenerate, then $X$ and $X_{\Delta}^{U}$ are topologically equivalent near infinity (see Section \ref{sec-upper} for precise definitions and a statement of this result). When $X_{\Delta}^{U}$ is Newton degenerate, an analogous result was proved in \cite{OliVal26} for a class of polynomial vector fields.

At this point, it is important to recall the following result. For a fixed Newton polytope $\mathcal{P}$, denote the set of all polynomial vector fields with upper principal part associated with $\mathcal{P}$ by $\mathfrak{U}(\mathcal{P})$, and the subset of all vector fields with Newton non-degenerate upper principal part by $\mathfrak{U}_{0}(\mathcal{P})$. Then we have the disjoint union $\mathfrak{U}(\mathcal{P}) = \mathfrak{U}_{0}(\mathcal{P}) \cup \mathfrak{U}_{1}(\mathcal{P})$, where $\mathfrak{U}_{1}(\mathcal{P})$ is the set of all polynomial vector fields with Newton degenerate upper principal part, and $\mathfrak{U}_{0}(\mathcal{P})$ is open and dense (see Proposition \ref{prop-generic-upper} for a precise statement).

This paper is a natural continuation of \cite{DOP24, OliVal26}. Our goal is to consider vector fields in $\mathfrak{U}_{1}(\mathcal{P})$, that is, planar polynomial vector fields whose upper principal part is Newton degenerate, and we generalize the results of \cite{DOP24, OliVal26}. We introduce the notion of the \textit{generalized upper principal part} of $X$, which will be denoted by $X_{\Gamma}^{U}$. Roughly speaking, $X_{\Gamma}^{U}$ is given by the upper principal part $X_{\Delta}^{U}$ plus some monomials related to points below the upper diagram of the Newton polytope (see Section \ref{sec-def-gupp} for a precise definition). Actually, one can do even better. Indeed, it is possible to drop some monomials of $X_{\Gamma}^{U}$, which leads us to the definition of the \textit{minimal generalized upper principal part} of $X$, denoted by $X_{G}^{U}$ and defined in Section \ref{sec-def-gupp}. We highlight that, when $X_{\Delta}^{U}$ is Newton non-degenerate, then the minimal generalized upper principal part $X_{G}^{U}$ coincides with the upper principal part $X_{\Delta}^{U}$.

Our main result, which is Theorem \ref{mthm} stated in Section \ref{sec-main-thms}, says that there exists a set $\widetilde{\mathfrak{U}}_{1}$, which is open and dense in $\mathfrak{U}_{1}$, such that, if $X \in \widetilde{\mathfrak{U}}_{1}$, then $X$ and its minimal generalized upper principal part $X_{G}^{U}$ are topologically equivalent near infinity, provided that the minimal generalized principal part $X_{G}^{U}$ satisfies some additional non-degeneracy hypothesis (see Definition \ref{def-mgupp-nd}). Therefore, we have the disjoint union $\mathfrak{U}_{1}(\mathcal{P}) = \widetilde{\mathfrak{U}}_{1}(\mathcal{P}) \cup \mathfrak{U}_{2}(\mathcal{P})$, and (in general) one cannot guarantee that a polynomial vector field $X \in \mathfrak{U}_{2}(\mathcal{P})$ is topologically equivalent to its minimal generalized upper principal part $X_{G}^{U}$. The proof of Theorem \ref{mthm} is given in Section \ref{sec-proofs}, and an alternative proof for the case where the upper diagram of the Newton polytope has only one \textit{useful segment} (see Section \ref{sec-upper} for the definition) is given in Section \ref{sec-one-useful-segment}. In both proofs, the sets $\widetilde{\mathfrak{U}}_{1}(\mathcal{P})$ and $\mathfrak{U}_{2}(\mathcal{P})$ are precisely defined, and they depend on the fixed Newton polygon $\mathcal{P}$ but not on the given vector field $X$.

It is important to highlight some features of the proof given in this paper. This paper is a natural continuation of \cite{DOP24, OliVal26}, and it is also motivated by \cite{Zup96,Zup00}, where the author performed a study near the origin and proved the topological equivalence of a real planar analytic vector field and its minimal generalized principal part $X_{G}^{L}$ (which will be called the \textit{minimal generalized lower principal part} in this paper, see Section \ref{sec-lower} for a precise definition). The idea of the proof given in \cite{Zup96,Zup00} relies on an inductive process based on the algorithm of resolution of singularities provided by Pelletier in \cite{Pel95} and on applying the Normal Form Theorem (see Appendix \ref{appendix-nft} and the references therein) in order to study the singularities along the exceptional divisor. In the proof given in the present paper, it is not necessary to carry out a strategy based on mathematical induction. Indeed, our strategy is to perform \textit{toric compactification} (see Section \ref{sec-proofs}) and then use the Normal Form Theorem as stated in Appendix \ref{appendix-nft} in order to study the singularities at infinity. This approach suggests an alternative proof for the main result of \cite{Zup96,Zup00} (see the concluding remarks in Section \ref{sec-concluding-remarks}).

This paper is structured as follows. The basic definitions that will be used throughout this paper are given in Section \ref{sec-main-definitions}. In Section \ref{sec-def-gupp}, we introduce the notions of Newton decomposition, generalized upper principal part, and minimal generalized upper principal part, and we provide some examples as well. The main result of the paper, Theorem \ref{mthm}, is stated in Section \ref{sec-main-thms}, and its proof is given in Section \ref{sec-proofs}. An alternative proof of Theorem \ref{mthm} is given in Section \ref{sec-one-useful-segment} for the particular case where the upper diagram has only one useful segment. Finally, in Appendices \ref{appendix-nft} and \ref{appendix-proofs}, we recall the Normal Form Theorem and prove some technical lemmas used in the proofs.

\section{Main definitions and tools}\label{sec-main-definitions}

Denote the set of all real planar polynomial vector fields by $\mathfrak{X}(\mathbb{R}^{2})$. We write $X\in\mathfrak{X}(\mathbb{R}^{2})$ as a sum of quasi-homogeneous polynomial vector fields in the so called \textit{logarithmic basis} (see \cite{Pan06})
\begin{equation}\label{eq-quasi-homogeneous-vf}
\begin{array}{rcl}
X(x,y) & = & \displaystyle\sum_{d=-1}^{\delta - 1}X^{(\alpha,\beta)}_{d}(x,y), \\
X^{(\alpha,\beta)}_{d}(x,y) & = & \displaystyle\sum_{\alpha m+\beta n=d}\left( a_{m,n}x^{m+1}y^{n}\displaystyle\frac{\partial}{\partial x} + b_{m,n}x^{m}y^{n+1}\displaystyle\frac{\partial}{\partial y} \right),
\end{array}
\end{equation}
with $a_{m,n},b_{m,n}\in\mathbb{R}$ and $m,n\in \mathbb{Z}$ satisfying: If $m < -1$ or $n \leq -1$, then $a_{m,n} = 0$; and if $m \leq -1$ or $n < -1$, then $b_{m,n} = 0$. We say that equation \eqref{eq-quasi-homogeneous-vf} is a \textit{$(\alpha,\beta)$-decomposition} of the planar vector field $X$, and the vector field $X^{(\alpha,\beta)}_{d}$
is called \textit{$d$-level} of the $(\alpha,\beta)$-decomposition of $X$.

Suppose that the origin $0\in\mathbb{R}^{2}$ is a singularity of $X$. We say that $0$ is a \textit{semi-hyperbolic singularity} if the linearization of $X$ at $0$ has only one eigenvalue with nonzero real part. If $0\in\mathbb{R}^{2}$ is hyperbolic or semi-hyperbolic singularity, then $0$ will be called \textit{elementary singularity}.

\subsection{Poincar\'e--Lyapunov compactification}

Let $\omega = (\alpha,\beta)$ be a \emph{weight vector} of positive integers satisfying $\operatorname{gcd}(\alpha,\beta) = 1$, and denote by $\operatorname{Cs}\theta$ and $\operatorname{Sn}\theta$ the unique solutions of the Cauchy problem
\begin{equation*}\label{eq-cauchy-problem}
\displaystyle\frac{d}{d\theta}\operatorname{Cs}\theta = -\operatorname{Sn}^{2\alpha - 1}\theta, \quad \displaystyle\frac{d}{d\theta}\operatorname{Sn}\theta = \operatorname{Cs}^{2\beta - 1}\theta, \qquad \operatorname{Cs}0 = 1, \operatorname{Sn}0 = 0;
\end{equation*}
which are analytic and periodic, and whose period $T$ is given by the formula
\begin{equation*}\label{eq-period-lyapunov}
T = \displaystyle\frac{2\alpha^{\frac{1-2\alpha}{2\alpha}}}{\beta^{\frac{1}{2\alpha}}} \int_{0}^{1}(1-t)^{\frac{1-2\alpha}{2\alpha}}t^{\frac{1-2\beta}{2\beta}}dt.
\end{equation*}

Such pair of functions also satisfy the equation $$\beta\operatorname{Sn}^{2\alpha}\theta + \alpha\operatorname{Cs}^{2\beta}\theta = \alpha.$$

Given a planar vector field $X$, the analytic vector field $\overline{X}$ obtained after the change of coordinates
\begin{equation}\label{eq-cs-sn}
x = \displaystyle\frac{\operatorname{Cs}\theta}{r^{\alpha}}, \ \quad \ y = \displaystyle\frac{\operatorname{Sn}\theta}{r^{\beta}}, \ \quad \ \theta\in\mathbb{S}^{1}, r\geq 0;
\end{equation}
and by a multiplication of a suitable power of $r$ is called \emph{Poincaré--Lyapunov compactification of $X$} (PLC for short) and it is defined in the \emph{Poincaré--Lyapunov disk $\mathbb{D}_{(\alpha,\beta)}$} (PL-disk for short). In this new phase space, the set $\{r = 0\}$ plays the role of infinity. The change of coordinates given in \eqref{eq-cs-sn} is called \textit{quasi-polar coordinates.} Usually, instead of dealing with quasi-polar coordinates, one studies the infinity using directional charts given by the change of coordinates \begin{multicols}{2} \noindent\begin{equation}\label{eq-plc-x-positive}
x = v^{-\alpha}, \ y = uv^{-\beta},
\end{equation}
\begin{equation}\label{eq-plc-x-negative}
x = -v^{-\alpha}, \ y = uv^{-\beta},
\end{equation}
\begin{equation}\label{eq-plc-y-positive}
x = uv^{-\alpha}, \ y = v^{-\beta},
\end{equation}
\begin{equation}\label{eq-plc-y-negative}
x = uv^{-\alpha}, \ y = -v^{-\beta}.
\end{equation}
\end{multicols}

Observe that in each equation above the variables $(u,v)$ have different meanings. Moreover, after such transformations, it is necessary to multiply the vector field by a suitable power of $v$. The Equations \eqref{eq-plc-x-positive} and \eqref{eq-plc-x-negative} are called \textit{compactification in the positive} and \textit{negative $x$-directions}, respectively. Similarly, Equations \eqref{eq-plc-y-positive} and \eqref{eq-plc-y-negative} are called \textit{compactification in the positive} and \textit{negative $y$-directions}, respectively. Moreover, in the new phase plane, the set $\{v = 0\}$ plays the role of infinity. The vector field obtained after performing on $X$ the change of coordinates \eqref{eq-plc-x-positive}, \eqref{eq-plc-x-negative}, \eqref{eq-plc-y-positive} and \eqref{eq-plc-y-negative} (and multiplication by a suitable power of $v$) will be denoted by $\overline{X}_{x}^{+}$, $\overline{X}_{x}^{-}$, $\overline{X}_{y}^{+}$ and $\overline{X}_{y}^{-}$, respectively. We refer to \cite[Chapters 5 and 9]{DLA06} for more details about the Poincaré--Lyapunov compactification.

The weight vector above can be chosen using the Newton polytope, which we define in the next Section.

\subsection{Newton polytope of polynomial planar vector fields}

Given the vector field \eqref{eq-quasi-homogeneous-vf}, we associate the monomials $a_{m,n}x^{m+1}y^{n}$ and $b_{m,n}x^{m}y^{n+1}$ with nonzero coefficients to the point $(m,n)$ in the plane of powers. Observe that each point $(m,n)$ is contained in a line of the form $\{(m,n); \alpha m + \beta n = d\}$, with $\omega = (\alpha,\beta)\in\mathbb{Z}^{2}$ satisfying $\gcd(\alpha,\beta) = 1$ and $d\in\mathbb{Z}$.

The \textit{support $\mathcal{S}_{X}$ of $X$} is the set
$$\mathcal{S}_{X} = \{(m,n)\in \mathbb{Z}^{2}; a_{m,n}^{2} + b_{m,n}^{2} \neq 0\},$$
which is finite because $X$ is polynomial. The \textit{Newton polytope $\mathcal{P}_{X}$} associated to the polynomial vector field $X$ is the convex hull of the support $\mathcal{S}_{X}$. The Newton polytope strongly depends on the coordinate system adopted. Moreover, when the vector field $X$ considered is well understood, for simplicity sake we denote $\mathcal{S}_{X}$ by $\mathcal{S}$ and $\mathcal{P}_{X}$ by $\mathcal{P}$. Since $X$ is a polynomial vector field, if $\mathcal{P}$ has empty interior, then $\mathcal{P}$ is either a point or a compact segment.

\subsection{Lower diagram and lower principal part}\label{sec-lower}

The \textit{lower diagram of $\mathcal{P}$} will be denoted as $\mathcal{P}^{L}$ and it is the union of all compact segments of the boundary of the set
$$\operatorname{conv}\left(\mathcal{P} + \left(\mathbb{R}_{\geq 0}\right)^{2}\right),$$
in which $\operatorname{conv}(\cdot)$ is the operation of convex closure, $\left(\mathbb{R}_{\geq 0}\right)^{2}$ denotes the closure of the first quadrant of $\mathbb{R}^{2}$ and $``+"$ is the Minkowski sum of convex polyhedrons. In many references, the lower diagram is also called \textit{Newton diagram}.

%The non-smooth points of the lower diagram $\mathcal{P}^{L}$ will be called \textit{vertices}. They are enumerated in the counterclockwise sense as $p_{0}^{L},\dots,p_{k}^{L}\in\mathcal{P}^{L}$ so that the vertex $p_{0}^{L}$ is the first point of the support $\mathcal{S}$ according to the lexicographical order.

The non-smooth points of $\mathcal{P}^{L}$ are called \textit{vertices} of $\mathcal{P}^{L}$. We say that a compact segment contained in $\mathcal{P}^{L}$ is \textit{useful} if it is not entirely contained in $[-1,0]\times \mathbb{Z}$ or $\mathbb{Z}\times [-1,0]$. Equivalently, a compact segment contained in $\mathcal{P}^L$ is \textit{useful} if it intersects the interior of the first quadrant of the plane of powers. The lower diagram $\mathcal{P}^{L}$ is the union of a finite number of compact segments, and they will be enumerated in the counterclockwise sense $\gamma_{0}^{L}, \dots, \gamma_{k+1}^{L}$ in such a way that, for $j = 0$ and $j = k+1$, the segment $\gamma_{j}^{L}$ is not useful (when it exists). In other words, the useful segments will be denoted as $\gamma_{1}^{L}, ..., \gamma_{k}^{L}$. See Figure \ref{fig-def-polytope}.

The \textit{lower principal part of $X$} is the polynomial vector field
$$X_{\Delta}^{L}(x,y) = \sum_{(m,n)\in\mathcal{P}^{L}} \left( a_{m,n}x^{m+1}y^{n}\displaystyle\frac{\partial}{\partial x} + b_{m,n}x^{m}y^{n+1}\displaystyle\frac{\partial}{\partial y} \right),$$
and we say that the lower principal part $X_{\Delta}^{L}$ is \textit{Newton non-degenerate} if any quasi-homogeneous component 
$$X_{\gamma_{j}^{L}}(x,y) = \sum_{(m,n)\in\gamma_{j}^{L}} \left( a_{m,n}x^{m+1}y^{n}\displaystyle\frac{\partial}{\partial x} + b_{m,n}x^{m}y^{n+1}\displaystyle\frac{\partial}{\partial y} \right),$$
has no singularities in $(\mathbb{R}^{*})^{2}$, in which $\gamma_{j}^{L}\subset\mathcal{P}^{L}$ (for $j = 0,\dots, k+1$) and $\mathbb{R}^{*} = \mathbb{R}\backslash\{0\}$. In other words, $X_{\Delta}^{L}$ is Newton non-degenerate if, for every $\gamma_{j}^{L}\subset\mathcal{P}^{L}$, the vector field $X_{\gamma_{j}^{L}}$ does not have singularities outside the coordinate axes. Otherwise, we say that $X_{\Delta}^{L}$ is \textit{Newton degenerate}. The following result is well known in the literature and it holds for analytic vector fields, however, we state it for polynomial vector fields, which is the scope of this paper.

\begin{proposition}[Berezovskaya \cite{Ber78}, Brunella and Miari \cite{BM90}]\label{prop-generic-lower}
Let $\mathcal{P}$ be a fixed Newton polytope. Denote the set of all polynomial vector fields with lower principal part associated to $\mathcal{P}$ by $\mathfrak{L}(\mathcal{P})$, and the subset of all polynomial vector fields with Newton non-degenerate lower principal part by $\mathfrak{L}_{0}(\mathcal{P})$. Then $\mathfrak{L}_{0}(\mathcal{P})$ is an open and dense subset of $\mathfrak{L}(\mathcal{P})$.
\end{proposition}

Denote the set of all polynomial vector fields with Newton degenerate lower principal part by $\mathfrak{L}_{1}(\mathcal{P})$. It follows from Proposition \ref{prop-generic-lower} that we have the disjoint union $\mathfrak{L}(\mathcal{P}) = \mathfrak{L}_{0}(\mathcal{P})\cup\mathfrak{L}_{1}(\mathcal{P})$, where $\mathfrak{L}_{0}(\mathcal{P})$ is open and dense (in fact, Zariski-open).

In \cite{Zup96,Zup00}, the author studied vector fields in the set $\mathfrak{L}_{1}(\mathcal{P})$ using the notion of \textit{generalized lower principal part} (which was called simply as generalized principal part in such references). \v{Z}upanovi\'c proved that there is an open and dense subset $\widetilde{\mathfrak{L}}_{1}(\mathcal{P})\subset \mathfrak{L}_{1}(\mathcal{P})$ such that, if the vector field $X$ belongs to $\widetilde{\mathfrak{L}}_{1}(\mathcal{P})$ (among other non degeneracy assumptions), then the $X$ and its generalized lower principal part are topologically equivalent. Therefore we have the disjoint union $\mathfrak{L}_{1}(\mathcal{P}) = \widetilde{\mathfrak{L}}_{1}(\mathcal{P})\cup \mathfrak{L}_{2}(\mathcal{P})$, in which (in general) it cannot be assured that a vector field in $\mathfrak{L}_{2}(\mathcal{P})$ is topologically equivalent to its generalized lower principal part.

\subsection{Upper diagram and upper principal part}\label{sec-upper}

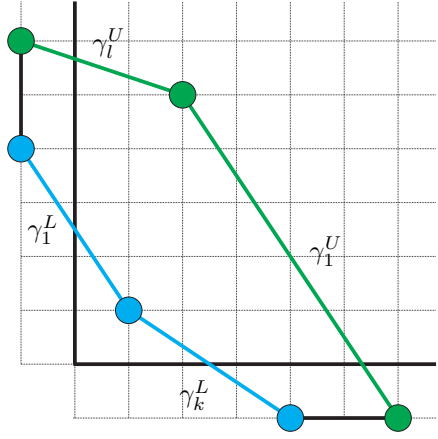
\begin{figure}[ht]\center{
\begin{overpic}[width=0.35\textwidth]{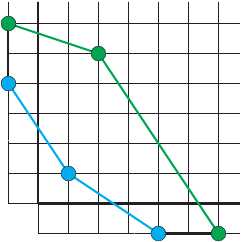}
%\begin{overpic}[grid,tics=10,width=0.35\textwidth]{fig-def-polytope.pdf}

\put(70,40){\footnotesize{$\gamma_{1}^{U}$}}
\put(5,45){\footnotesize{$\gamma_{1}^{L}$}}

\put(20,88){\footnotesize{$\gamma_{l}^{U}$}}
\put(40,7){\footnotesize{$\gamma_{k}^{L}$}}

\end{overpic}}
\caption{\footnotesize{Lower and upper boundaries $\mathcal{P}^{L}$ (highlighted in blue) and $\mathcal{P}^{U}$ (highlighted in green), respectively.}}
\label{fig-def-polytope}
\end{figure}

One can state similar definitions in order to study the dynamics at infinity, which is the scope of this paper. This lead us to the notion of upper diagram and upper principal part. The \textit{upper diagram of $\mathcal{P}$} will be denoted as $\mathcal{P}^{U}$ and it is the union of all compact segments of the boundary of the set
$$\operatorname{conv}\left(\mathcal{P} \cup \{(0,-1);(-1,0)\}\right)$$
that do not contain $(0,-1)$ nor $(-1,0)$. Just as in Section \ref{sec-lower}, the non-smooth points of $\mathcal{P}^{U}$ are called \textit{vertices} of $\mathcal{P}^{U}$. The notion of useful segment is defined in a completely analogous way as in Section \ref{sec-lower}. The upper diagram $\mathcal{P}^{U}$ is the union of a finite number of compact segments, and they will also be enumerated in the counterclockwise sense $\gamma_{0}^{U}, \dots, \gamma_{l+1}^{U}$ in such a way that, for $j = 0$ or $j = l+1$, the segment $\gamma_{j}^{U}$ is not useful (when it exists). Therefore, the useful segments are $\gamma_{1}^{U},\dots,\gamma_{l}^{U}$. See Figure \ref{fig-def-polytope}.

\begin{remark}
{\rm The definition of upper diagram presented in the last paragraph agrees with \cite[Definition 1]{DOP24}.}    
\end{remark}

Similarly to Section \ref{sec-lower}, it is possible to define the notion of \textit{upper principal part of $X$} (see \cite{DOP24}), which is the polynomial vector field
\begin{equation}\label{eq-upper-pp}
X_{\Delta}^{U}(x,y) = \sum_{(m,n)\in\mathcal{P}^{U}} \left( a_{m,n}x^{m+1}y^{n}\displaystyle\frac{\partial}{\partial x} + b_{m,n}x^{m}y^{n+1}\displaystyle\frac{\partial}{\partial y} \right).
\end{equation}

The upper principal part $X_{\Delta}^{U}$ is \textit{Newton non-degenerate} if any quasi-homogeneous component 
\begin{equation}\label{eq-quasi-upper}
X_{\gamma_{j}^{U}}(x,y) = \sum_{(m,n)\in\gamma_{j}^{U}} \left( a_{m,n}x^{m+1}y^{n}\displaystyle\frac{\partial}{\partial x} + b_{m,n}x^{m}y^{n+1}\displaystyle\frac{\partial}{\partial y} \right),
\end{equation}
has no singularities in $(\mathbb{R}^{*})^{2}$, in which $\gamma_{j}^{U}\subset\mathcal{P}^{U}$, for $j = 0, \dots, l+1$. One can prove an analogous result to Proposition \ref{prop-generic-lower}.

\begin{proposition}[see \cite{DOP24}]\label{prop-generic-upper}
Let $\mathcal{P}$ be a fixed Newton polytope. Denote the set of all polynomial vector fields with upper principal part associated to $\mathcal{P}$ by $\mathfrak{U}(\mathcal{P})$, and the subset of all polynomial vector fields with Newton non-degenerate upper principal part by $\mathfrak{U}_{0}(\mathcal{P})$. Then the set $\mathfrak{U}_{0}(\mathcal{P})$ is an open and dense subset of $\mathfrak{U}(\mathcal{P})$.
\end{proposition}

It follows from Proposition \ref{prop-generic-upper} we have the disjoint union $\mathfrak{U}(\mathcal{P}) = \mathfrak{U}_{0}(\mathcal{P})\cup\mathfrak{U}_{1}(\mathcal{P})$, where $\mathfrak{U}_{1}(\mathcal{P})$ is the set of all polynomial vector fields with Newton degenerate upper principal part and $\mathfrak{U}_{0}(\mathcal{P})$ is an open and dense subset of $\mathfrak{U}(\mathcal{P})$.

The goal of this paper is to study vector fields in the set $\mathfrak{U}_{1}(\mathcal{P})$ using the notion of \textit{generalized upper principal part}. More precisely, we will prove that there is an open and dense subset $\widetilde{\mathfrak{U}}_{1}(\mathcal{P})\subset \mathfrak{U}_{1}(\mathcal{P})$ such that, if the polynomial vector field $X$ belongs to this set and it satisfies some additional non degeneracy conditions, then $X$ is topologically equivalent to its generalized upper principal part. Therefore, we have the disjoint union $\mathfrak{U}_{1}(\mathcal{P}) = \widetilde{\mathfrak{U}}_{1}(\mathcal{P})\cup \mathfrak{U}_{2}(\mathcal{P})$, and (in general) one cannot assure that a vector field on $\mathfrak{U}_{2}(\mathcal{P})$ is topologically equivalent to its generalized upper principal part. It is also important to remark that the subsets $\widetilde{\mathfrak{U}}_{1}(\mathcal{P})$ and $\mathfrak{U}_{2}(\mathcal{P})$ of $\mathfrak{U}_{1}(\mathcal{P})$ will not depend on the initial vector field $X$.

The main result of this paper is Theorem \ref{mthm}, which is stated in Section \ref{sec-main-thms} and proved in Section \ref{sec-proofs}. Firstly, we introduce some notation and definitions in Section \ref{sec-def-gupp}.

\section{Newton decomposition and generalized upper principal part}\label{sec-def-gupp}

This section is devoted to define the notion of Newton decomposition associated to $\mathcal{P}^{U}$ (in which $\mathcal{P}$ is a fixed Newton polytope), as well as the notion of generalized upper principal part, which takes into account vector fields in a subset of $\mathfrak{U}_{1}(\mathcal{P})$.

Firstly, let us introduce the notion of \textit{Newton decomposition}. Consider a planar polynomial vector field $X$ and denote its Newton polytope by $\mathcal{P}$. It is possible to decompose $X$ as a sum of polynomial vector fields in a slightly different way than \eqref{eq-quasi-homogeneous-vf}. In what follows, it is important to keep in mind the notation established in Section \ref{sec-upper}, and the following construction is illustrated by the examples in Sections \ref{sec-example} and \ref{sec-another-ex}.

Let $\gamma_{j}^{U}\subset \mathcal{P}^{U}$ be a compact segment, for $j = 0,\dots, l+1$. Observe that each $\gamma_{j}^{U}\subset \mathcal{P}^{U}$ is contained in a line of the form $\{(m,n); \alpha_{j} m + \beta_{j} n = \delta_{j}\}$ (for $j = 0,\dots, l+1$), and we suppose that, in such an expression of the support line, the vector $\omega_{j} = (\alpha_{j},\beta_{j})\in\mathbb{Z}^{2}$ satisfies 
$$\operatorname{gcd}(\alpha_{j},\beta_{j}) = 1, \quad \text{and} \quad \delta_{j} = \displaystyle\min_{(m,n)\in\mathcal{S}}\{\alpha_{j} m + \beta_{j} n \} \leq 0.$$

The vector $\omega_{j} = (\alpha_{j},\beta_{j})$ satisfying such conditions will be called \textit{inward normal vector}, and the collection of all inward normal vectors associated to the upper diagram $\mathcal{P}^{U}$ is denoted by $\Sigma_{\mathcal{P}}^{*} = \{\omega_{1},\dots,\omega_{l}\}$. 

For each inward normal vector $\omega_{j} = (\alpha_{j},\beta_{j})$, one can write
$$
\delta_{j}^{(0)} = \displaystyle\min_{(m,n)\in\mathcal{S}}\{\alpha_{j} m + \beta_{j} n \}, \ \quad \
   \Gamma_{j}^{(0)}  =  \{(m,n)\in\mathcal{S} \ ; \ \alpha_{j} m + \beta_{j} n = \delta_{j}^{(0)}\},
$$
for each $j = 0,\dots, l+1$, where $\mathcal{S}$ is the support of $X$.  It is clear that $\Gamma_{j}^{(0)} = \gamma_{j}^{U}\cap\mathcal{S}$. Moreover, the integer $\delta_{j}^{(0)} \leq 0$ will be called \textit{degree of quasi homogeneity} of $\Gamma_{j}^{(0)}$. Setting $\Gamma^{(0)} = \cup_{j = 0}^{l+1}\Gamma_{j}^{(0)}$, define the vector field
$$
X_{0}^{U}(x,y) = \sum_{(m,n)\in\Gamma^{(0)}} \left( a_{m,n}x^{m+1}y^{n}\displaystyle\frac{\partial}{\partial x} + b_{m,n}x^{m}y^{n+1}\displaystyle\frac{\partial}{\partial y} \right).
$$

Since $\Gamma^{(0)} = \mathcal{P}^{U}\cap\mathcal{S}$ it is straightforward to see that $X_{0}^{U} = X_{\Delta}^{U}$, that is, the vector field $X_{0}^{U}$ coincides with the upper principal part \eqref{eq-upper-pp}. Now, define
$$
   \delta_{j}^{(1)}  =  \displaystyle\min_{(m,n)\in\mathcal{S}\backslash \Gamma_{j}^{(0)}}\{\alpha_{j} m + \beta_{j} n \}, \ \quad \
   \Gamma_{j}^{(1)}  =  \{(m,n)\in\mathcal{S}\backslash \Gamma_{j}^{(0)} \ ; \ \alpha_{j} m + \beta_{j} n = \delta_{j}^{(1)}\},
$$
for all $j = 0,\dots, l+1$. It can be checked that $\delta_{j}^{(0)} \leq \delta_{j}^{(1)} \leq 0$. Geometrically, the points of $\Gamma_{j}^{(1)}$ are contained in the first non-empty parallel line below the segment $\gamma_{j}^{U}$. Moreover, the set  $\Gamma_{j}^{(1)}$ might contain points of $\Gamma^{(0)}$. Then, define the set
\begin{equation*}
\Gamma^{(1)} = \Big{(}\bigcup_{j = 0}^{l+1}\Gamma_{j}^{(1)}\Big{)}\backslash \Gamma^{(0)}.
\end{equation*}

The set $\Gamma^{(1)}$ is either empty, or it contains points of the support $\mathcal{S}$ that are below $\mathcal{P}^{U}$. We then define the vector field
$$X_{1}^{U}(x,y) = \sum_{(m,n)\in\Gamma^{(1)}} \left( a_{m,n}x^{m+1}y^{n}\displaystyle\frac{\partial}{\partial x} + b_{m,n}x^{m}y^{n+1}\displaystyle\frac{\partial}{\partial y} \right).$$

It can be checked that $X_{1}^{U}$ is either identically zero (when $\Gamma^{(1)}$ is empty), or its monomials are related to points below $\mathcal{P}^{U}$. By proceeding inductively with this construction, we define 
\begin{equation*}
   \delta_{j}^{(i)}  =  \displaystyle\min_{(m,n)\in\mathcal{S}\backslash \bigcup\Gamma_{j}^{(i-1)}}\{\alpha_{j} m + \beta_{j} n \}, \ \quad \
   \Gamma_{j}^{(i)}  =  \{(m,n)\in\mathcal{S}\backslash \cup\Gamma_{j}^{(i-1)} \ ; \ \alpha_{j} m + \beta_{j} n = \delta_{j}^{(i)}\},    
\end{equation*}
for each $j = 0,\dots, l+1$. Then, we define the set
\begin{equation*}
\Gamma^{(i)} = \Big{(}\displaystyle\bigcup_{j = 0}^{l+1}\Gamma_{j}^{(i)}\Big{)}\backslash \Big{(}\displaystyle\bigcup_{r = 0}^{i-1}\Gamma^{(r)}\Big{)}, 
\end{equation*}
and we finally define the vector field
$$X_{i}^{U}(x,y) = \sum_{(m,n)\in\Gamma^{(i)}} \left( a_{m,n}x^{m+1}y^{n}\displaystyle\frac{\partial}{\partial x} + b_{m,n}x^{m}y^{n+1}\displaystyle\frac{\partial}{\partial y} \right).$$

\begin{figure}[ht]\center{
\begin{overpic}[width=0.85\textwidth]{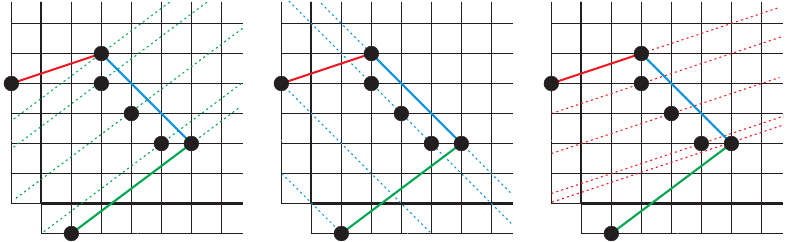}
%\begin{overpic}[grid,tics=10,width=0.85\textwidth]{fig-decomposition.pdf}
\put(12,1){\footnotesize{$\gamma_{j-1}^{U}$}}
\put(27,12){\footnotesize{$\Gamma_{j-1}^{(0)}$}}
\put(30,18){\footnotesize{$\Gamma_{j-1}^{(1)}$}}
\put(30,25){\footnotesize{$\Gamma_{j-1}^{(2)}$}}

\put(52,20){\footnotesize{$\gamma_{j}^{U}$}}
\put(62,10){\footnotesize{$\Gamma_{j}^{(0)}$}}
\put(57,4){\footnotesize{$\Gamma_{j}^{(1)}$}}
\put(50,0){\footnotesize{$\Gamma_{j}^{(2)}$}}

\put(74,24){\footnotesize{$\gamma_{j+1}^{U}$}}
\put(85,28){\footnotesize{$\Gamma_{j+1}^{(0)}$}}
\put(100,25){\footnotesize{$\Gamma_{j+1}^{(1)}$}}
\put(100,20){\footnotesize{$\Gamma_{j+1}^{(2)}$}}

\end{overpic}}
\caption{\footnotesize{Sketch of the sets $\Gamma^{(i)}_{j}$, which are contained in the dashed lines. The index $j$ refers to a compact segment of the upper diagram, and the index $i$ refers to the level of the set.}}
\label{fig-decomposition}
\end{figure}
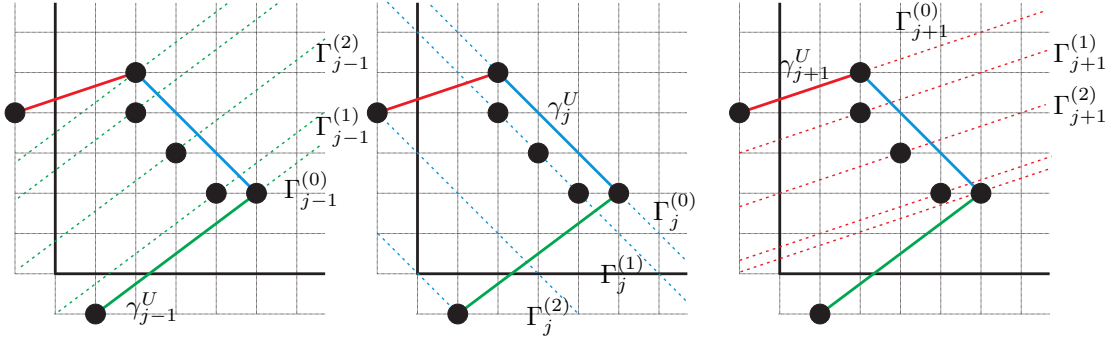

Observe that it is possible that a point $(m,n)\in\mathcal{S}$ may be contained in $\Gamma_{j_{1}}^{(i)}\cap\Gamma_{j_{2}}^{(i)}$, for $0\leq j_{1,2}\leq l+1$ with $j_{1} \neq j_{2}$ and some positive integer $i$. See Figure \ref{fig-decomposition}.

\begin{definition}
Let $N$ be the smallest positive integer such that $\cup_{i = 0}^{N}\Gamma^{(i)} = \mathcal{S}$. The \textit{Newton decomposition} of the planar polynomial vector field $X$ with respect to the upper diagram $\mathcal{P}^{U}$ is given by
\begin{equation*}
X(x,y) = \displaystyle\sum_{i = 0}^{N}X_{i}^{U}(x,y).    
\end{equation*}
The \textit{generalized upper principal part} of $X$ with respect to the upper diagram $\mathcal{P}^{U}$ is the planar polynomial vector field
\begin{equation*}
     X_{\Gamma}^{U} = X_{0}^{U} + X_{1}^{U}  = X_{\Delta}^{U} + X_{1}^{U}.  
\end{equation*}
\end{definition}

\begin{remark}
We refer to \cite[Section 2]{Zup00} for the definition of generalized lower principal part (which was called simply as generalized principal part in such reference). It is important that the reader do not confuse the Newton decomposition as a vector field obtained from the so called \textit{Newton filtration} as defined in \cite[Section 2.1]{BivHua19} (see also \cite[Section 2]{Kou76}).
\end{remark}

The goal of this paper is to prove that, under non degeneracy conditions, the vector fields $X$ and $X^{U}_{\Gamma}$ are topologically equivalent in a neighborhood of infinity. However, there are some monomials of $X^{U}_{\Gamma}$ that can be dropped. Therefore, in what follows, our goal is to define the notion of \textit{minimal} generalized upper principal part.

%In what follows, the \textit{cone over $\gamma_{j}^{U}$} is defined as the set $C(\gamma_{j}^{U}) = \operatorname{conv}(\gamma_{j}^{U}\cup\{0\})$. 

\begin{definition}\label{def-mgupp}
Let $X$ be a planar polynomial vector field and $\mathcal{P}$ its Newton polygon. The \textit{minimal generalized upper principal part} of $X$ is the vector field $X_{G}^{U}$ obtained by dropping some monomials from $X_{1}^{U}$ of the generalized upper principal part $X_{\Gamma}^{U}$ as follows. Consider $(m,n)\in \Gamma^{(1)}$. The monomials of $X_{\Gamma}^{U}$ related to $(m,n)\in\Gamma^{(1)}$ are dropped if, for every $\Gamma_{j}^{(1)}$ containing $(m,n)$, the vector field $X_{\gamma_{j}^{U}}$ defined in \eqref{eq-quasi-upper} does not have singularities in $(\mathbb{R}^{*})^{2}$.
\end{definition}

Intuitively speaking, one can interpret Definition \ref{def-mgupp} as follows. If the vector field $X_{\gamma_{j}^{U}}$ defined in \eqref{eq-quasi-upper} is degenerate in the sense that it has singularities in $(\mathbb{R}^{*})^{2}$, then $X_{G}^{U}$ must have all the monomials associated to the points of the set $\Gamma_{j}^{(1)}$. Concerning the non useful segments $\gamma_{0}^{U}$ and $\gamma_{l+1}^{U}$, observe that the vector fields $X_{\gamma_{0}^{U}}$ and $X_{\gamma_{l+1}^{U}}$ do not have singularities in $(\mathbb{R}^{*})^{2}$. Therefore, monomials related to points of $\Gamma_{0}^{(1)}$ and $\Gamma_{l+1}^{(1)}$ that do not belong to any other $\Gamma_{j}^{(1)}$ for $j = 1,\dots,l$ can be dropped.

\begin{remark}
It can be easily checked that, if $\gamma_{j}^{U}$ is not a useful segment (for $j = 0$ or $j = l+1$), then the singularities of $X_{\gamma_{j}^{U}}$ are not contained in $(\mathbb{R}^{*})^{2}$. Therefore, it suffices to verify whether $X_{\Delta}^{U}$ is Newton degenerate or not by looking to the useful segment.
\end{remark}

Finally, it is also important to remark that, if $X_{\gamma_{j}^{U}}$ has no singularities in $(\mathbb{R}^{*})^{2}$ for every compact segment $\gamma_{j}^{U}\subset\mathcal{P}^{U}$ (that is, $X_{\Delta}^{U}$ is Newton non-degenerate), then it follows directly from Definition \ref{def-mgupp} that the minimal generalized upper principal part $X_{G}^{U}$ coincides with the upper principal part $X_{\Delta}^{U}$.

% Suppose that $\mathcal{P}$ has non-empty interior. Observe that $X_{G}^{U}$ might have monomials related to points in the interior of $\mathcal{P}$ that are in the diagram of the cone $C(\gamma_{j}^{U})$, for some $\gamma_{j}^{U}\subset\mathcal{P}^{U}$. 

%We follow the presentation of \cite{BivHua19}. 

%Denote $M = \operatorname{lcm}\{|\delta_{1}|,\dots,|\delta_{l}|\}$, then $M > 0$. For each $j = 1,\dots,l$, define the linear map $\phi_{j}:\mathbb{R}^{2}\rightarrow \mathbb{R}$ as
%$$\phi_{j}(p) = \frac{M}{\delta_{j}}\langle \omega_{j}, p \rangle.$$

%The \textit{filtrating map associated to $\mathcal{P}^{U}$} is the map $\phi:\mathbb{R}^{2}\rightarrow \mathbb{R}$ given by 
%$$\phi(p) = \displaystyle\operatorname{max} \{ \phi_{j}(p); 1 \leq j \leq l\}.$$

%\begin{lemma}\cite[Lemma 2.2]{BivHua19}
%The filtrating map $\phi$ associated to $\mathcal{P}^{U}$ satisfies the following properties.
%\begin{itemize}
%    \item[(a)] $\phi(\mathbb{Z}_{\geq 0}^{2}) \subset \mathbb{Z}_{\geq 0}.$
%    \item[(b)] Let $p\in\mathbb{R}^{2}$. Then $\phi(p) = \phi_{j}(p)$ if, and only if, $p\in\ C(\gamma_{j}^{U})$ for some $j = 1, \dots, l$.
%    \item[(c)] The map $\phi$ is linear in each cone $C(\gamma_{j}^{U})$, in which $\gamma_{j}^{U}\in\mathcal{P}^{U}$, $j = 1,\dots, l$. 
%    \item[(d)] Given $p,q\in\mathbb{R}^{2}$, then $\phi(p + q) \leq \phi(p) + \phi(q)$. Moreover, the equality holds if, and only if, $p$ and $q$ belong to the same cone $C(\gamma_{j}^{U})$, for some $j = 1,\dots, l$.
%\end{itemize}
%\end{lemma}

\subsection{An example}\label{sec-example}

Before we state the main Theorems of this paper in Section \ref{sec-main-thms}, we consider an example to illustrate all the definitions so far.

Consider the polynomial vector field
\begin{equation}\label{eq-example}
\begin{array}{rcl}
    X(x,y) & = & \big{(}ay^{5} + b_{1}x^{4}y^{3} + d_{1}x^{4}y^{2} + e_{1}x^{2}y^{3} + f_{1}x^{3}y + gy^{4}\big{)}\displaystyle\frac{\partial}{\partial x}  \\
     & + &\big{(}cx^{5} + b_{2}x^{3}y^{4}+d_{2}x^{3}y^{3}+e_{2}xy^{4}  + f_{2}x^{2}y^{2}\big{)}\displaystyle\frac{\partial}{\partial y}. 
\end{array}
\end{equation}

Suppose that the coefficients $a,c,d_{i},e_{i},f_{i},g$ of the vector field \eqref{eq-example} are all non zero, for $i = 1,2$. Moreover, assume $a\neq0$, $c\neq 0$ and $b_{1}^{2} + b_{2}^{2} \neq 0$. The support of $X$ is the set
$$\mathcal{S} = \{(-1,5), \ (5,-1), \ (3,3), \ (3,2), \ (1,3), \ (2,1), \ (-1,4)\}.$$

\begin{figure}[ht]\center{
\begin{overpic}[width=0.35\textwidth]{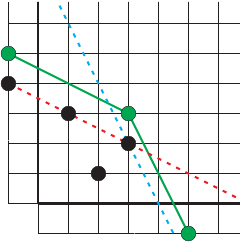}
%\begin{overpic}[grid,tics=10,width=0.35\textwidth]{fig-example.pdf}
\put(71,20){\footnotesize{$\gamma_{1}^{U}$}}
\put(58,5){\footnotesize{$\Gamma^{(1)}_{1}$}}

\put(25,70){\footnotesize{$\gamma_{2}^{U}$}}
\put(90,26){\footnotesize{$\Gamma^{(1)}_{2}$}}

\end{overpic}}
\caption{\footnotesize{Support $\mathcal{S}$ and upper diagram $\mathcal{P}^{U}$ (highlighted in green) of the vector field \eqref{eq-example}} in Section \ref{sec-example}. The points of $\Gamma_{1}^{(1)}$ belong to the first non empty parallel below $\gamma_{1}^{U}$, which is represented by the dashed blue line. On the other hand, points contained in $\Gamma_{2}^{(1)}$ belong to the first non empty parallel below $\gamma_{2}^{U}$, which is represented by the dashed red line. In this example, $\Gamma_{1}^{(1)}\cap\Gamma_{2}^{(1)} = \{(3,2)\}$.}
\label{fig-example}
\end{figure}
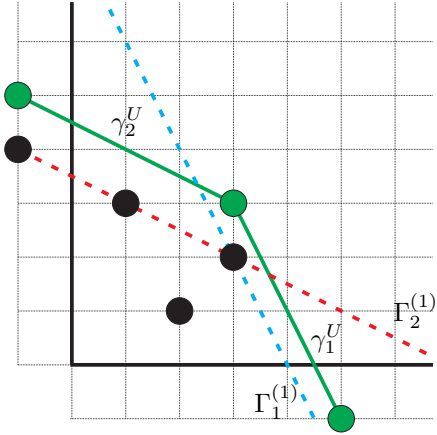

We refer to Figure \ref{fig-example}. The upper diagram $\mathcal{P}^{U}$ has two compact segments $\gamma_{1}^{U}$ and $\gamma_{2}^{U}$, and they are contained in the lines
\begin{equation*}
r_{1} = \{(m,n) \ ; \ -2m - n = -9\}, \quad r_{2} = \{(m,n) \ ; \ -m - 2n = -9\}. 
\end{equation*}

The inward normal vectors are $\omega_{1} = (-2,-1)$ and $\omega_{2} = (-1,-2)$. Moreover, it can be checked that $\delta_{1}^{(0)} = \delta_{2}^{(0)} = -9$, $\Gamma_{1}^{(0)} = \{(5,-1), (3,3)\}$ and $\Gamma_{2}^{(0)} = \{(-1,5), (3,3)\}$, which implies $\Gamma^{(0)} = \{(5,-1), (-1,5), (3,3)\}$. Therefore, one obtains the polynomial vector field
\begin{equation*}\label{eq-example-X0}
X_{0}^{U}(x,y) = X_{\Delta}^{U}(x,y) = \big{(}ay^{5} + b_{1}x^{4}y^{3}\big{)}\frac{\partial}{\partial x} +  \big{(}cx^{5} + b_{2}x^{3}y^{4}\big{)}\frac{\partial}{\partial y}.
\end{equation*}

Continuing this reasoning, one has that $\delta_{1}^{(1)} = -8$ and $\delta_{2}^{(1)} = -7$, and then $\Gamma_{1}^{(1)} = \{(3,2)\}$ and $\Gamma_{2}^{(1)} = \{(-1,4), (1,3), (3,2)\}$, and finally $\Gamma^{(1)} = \{(-1,4), (1,3), (3,2)\}$ (see Figure \ref{fig-example}). Therefore, one obtains the polynomial vector field
\begin{equation*}\label{eq-example-X1}
X_{1}^{U}(x,y) = \big{(}d_{1}x^{4}y^{2} + e_{1}x^{2}y^{3} + gy^{4}\big{)}\displaystyle\frac{\partial}{\partial x} + \big{(}d_{2}x^{3}y^{3}+e_{2}xy^{4}\big{)}\displaystyle\frac{\partial}{\partial y},
\end{equation*}
and finally we obtain the generalized upper principal part
\begin{equation*}\label{eq-example-xgen}
\begin{array}{rcl}
    X_{\Gamma}^{U}(x,y) & = & \big{(}ay^{5} + b_{1}x^{4}y^{3} + d_{1}x^{4}y^{2} + e_{1}x^{2}y^{3} + gy^{4}\big{)}\displaystyle\frac{\partial}{\partial x}  \\
     & + &\big{(}cx^{5} + b_{2}x^{3}y^{4}+d_{2}x^{3}y^{3}+e_{2}xy^{4} \big{)}\displaystyle\frac{\partial}{\partial y}. 
\end{array}
\end{equation*}

As discussed in Section \ref{sec-def-gupp}, one can further drop some monomials of $X_{\Gamma}^{U}$, and this will depend on the vector fields
\begin{equation*}
\begin{array}{rcrcl}
    X_{\gamma_{1}^{U}}(x,y) & = & b_{1}x^{4}y^{3}\frac{\partial}{\partial x} & + & \big{(}cx^{5} + b_{2}x^{3}y^{4}\big{)}\frac{\partial}{\partial y}, \\ X_{\gamma_{2}^{U}}(x,y) & = & \big{(}ay^{5} + b_{1}x^{4}y^{3}\big{)}\frac{\partial}{\partial x} & + & b_{2}x^{3}y^{4}\frac{\partial}{\partial y}.
\end{array}
\end{equation*}

If $b_{1} \neq0 \neq b_{2}$, then for $i = 1,2$ the vector field $X_{\gamma_{i}^{U}}$ does not have singularities in $(\mathbb{R}^{*})^{2}$, that is, the upper principal part $X_{\Delta}^{U}$ is Newton non-degenerate. Thus one can discard all the monomials associated with the sets $\Gamma_{i}^{(1)}$, which implies $X_{G}^{U} = X_{\Delta}^{U}$.

Now, assume that $b_{1} = 0$ and $b_{2} \neq 0$. In this case, $X_{\gamma_{1}^{U}}$ might have singularities in $(\mathbb{R}^{*})^{2}$, whereas $X_{\gamma_{2}^{U}}$ does not. Then one should keep only the monomial associated to $\Gamma_{1}^{(1)}$, and we drop the remaining monomials in the interior of $\mathcal{P}$. This gives us the minimal generalized upper principal part
\begin{equation*}
X_{G}^{U}(x,y) = \big{(}ay^{5} + b_{1}x^{4}y^{3} + d_{1}x^{4}y^{2}\big{)}\frac{\partial}{\partial x} +  \big{(}cx^{5} + b_{2}x^{3}y^{4} + d_{2}x^{3}y^{3}\big{)}\frac{\partial}{\partial y}.    
\end{equation*}

On the other hand, if $b_{1} \neq 0$ and $b_{2} = 0$, then $X_{\gamma_{2}^{U}}$ might have singularities in $(\mathbb{R}^{*})^{2}$, whereas $X_{\gamma_{1}^{U}}$ does not. In this case, we keep all the monomials associated to $\Gamma_{2}^{(1)}$, which gives us the minimal generalized upper principal part
\begin{equation*}
\begin{array}{rcl}
    X_{G}^{U}(x,y) = X_{\Gamma}^{U}(x,y) & = & \big{(}ay^{5} + b_{1}x^{4}y^{3} + d_{1}x^{4}y^{2} + e_{1}x^{2}y^{3} + gy^{4}\big{)}\displaystyle\frac{\partial}{\partial x}  \\
     & + &\big{(}cx^{5} + b_{2}x^{3}y^{4}+d_{2}x^{3}y^{3}+e_{2}xy^{4} \big{)}\displaystyle\frac{\partial}{\partial y}. 
\end{array}   
\end{equation*}

\subsection{Another example}\label{sec-another-ex}

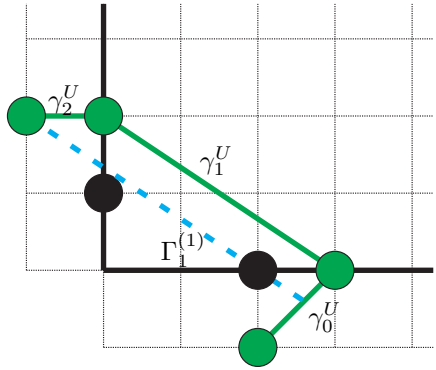
\begin{figure}[ht]\center{
\begin{overpic}[width=0.35\textwidth]{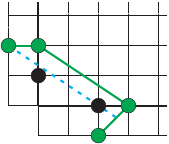}
%\begin{overpic}[grid,tics=10,width=0.35\textwidth]{fig-example-one-seg.pdf}
\put(70,10){\footnotesize{$\gamma_{0}^{U}$}}
\put(45,46){\footnotesize{$\gamma_{1}^{U}$}}
\put(10,60){\footnotesize{$\gamma_{2}^{U}$}}
\put(36,25){\footnotesize{$\Gamma^{(1)}_{1}$}}
\end{overpic}}
\caption{\footnotesize{Upper diagram $\mathcal{P}^{U}$ (highlighted in green) of the vector field $X$ given in Section \ref{sec-another-ex}. The support line of $\Gamma_{1}^{(1)} = \{(-1,2);(2,0)\}$ below the unique useful segment $\gamma_{1}^{U}$ is represented by a dashed blue line.}}
\label{fig-example-one-seg}
\end{figure}

Consider the vector field
\begin{equation*}\label{eq-ex-one-seg}
\begin{array}{rcl}
 \scalebox{0.97}{
 $X(x,y) = \left( ay^{2} + a_{1}xy^{2} + a_{2}x^{4} + a_{3}x^{3} + a_{4}xy\right)\displaystyle\frac{\partial}{\partial x} + \left(bx^{2} + b_{1}y^{3} + b_{2}x^{3}y + b_{3}x^{2}y + b_{4}y^{2}\right)\displaystyle\frac{\partial}{\partial y},$}
 \end{array}
\end{equation*}
whose support is $\mathcal{S} = \{(-1,2); \ (2,-1); \ (0,2); \ (3,0); \ (2,0); \ (0,1)\}$. The upper diagram $\mathcal{P}^{U}$ is the union of three compact segments, only one of them being useful, which is $\gamma_{1}^{U}$. This segment contains the points $(0,2)$ and $(3,0)$, and the support line of $\gamma_{1}^{U}$ is given by $r = \{(m,n); \ -2m - 3n = -6\}$. See Figure \ref{fig-example-one-seg}. 
Observe that 
\begin{equation*}
X_{\gamma_{1}^{U}}(x,y) = x\left( a_{1}y^{2} + a_{2}x^{3} \right)\displaystyle\frac{\partial}{\partial x} + y\left(b_{1}y^{2} + b_{2}x^{3}\right)\displaystyle\frac{\partial}{\partial y},    
\end{equation*}
and such a vector field admits singularities in $(\mathbb{R}^{*})^{2}$ if and only if $a_{1}b_{2} = a_{2}b_{1}$. If this is the case, then one should consider the minimal generalized upper principal part $X_{G}^{U}$. It can be checked that, following the notation introduced in Section \ref{sec-def-gupp}, $\Gamma_{1}^{(1)}$ contains $(-1,2)$ and $(2,0)$, whereas $\Gamma^{(1)}$ contains only $(2,0)$. The support line of $\Gamma_{1}^{(1)}$ is $\tilde{r} = \{(m,n); \ -2m -3n = -4\}$. Therefore, one obtains
\begin{equation*}
X_{G}^{U}(x,y) = \left( ay^{2} + a_{1}xy^{2} + a_{2}x^{4} + a_{3}x^{3} \right)\displaystyle\frac{\partial}{\partial x} + \left(bx^{2} + b_{1}y^{3} + b_{2}x^{3}y + b_{3}x^{2}y\right)\displaystyle\frac{\partial}{\partial y}.  
\end{equation*}

\section{Statement of the main results}\label{sec-main-thms}

This section is devoted to state the main theorem of this paper, which is Theorem \ref{mthm}. In what follows, denote the boundary of the Poincaré--Lyapunov disk $\mathbb{D}_{(\alpha,\beta)}$ by $\mathcal{D}$. With this notation, the set $\mathcal{D}$ plays the role of infinity of the phase plane. In addition, in the statement of Theorem \ref{mthm}, the subsets $\mathfrak{U}(\mathcal{P})$, $\mathfrak{U}_{0}(\mathcal{P})$ and $\mathfrak{U}_{1}(\mathcal{P})$ were defined in Section \ref{sec-upper}.

\begin{definition}[Definition 9, \cite{DOP24}]\label{def-free-charact}
A planar polynomial vector field $Y$ is \textit{free of characteristic orbits near infinity} if, for its PL-compactification $\overline{Y}$, there is no singularity in $\mathcal{D}$ having characteristic orbit intersecting the interior of $\mathbb{D}_{(\alpha,\beta)}$.     
\end{definition}

Geometrically, Definition \ref{def-free-charact} says that $Y$ is free of characteristic orbits near infinity if it presents a monodromic behavior near the boundary of $\mathbb{D}_{(\alpha,\beta)}$.

\begin{definition}\label{def-mgupp-nd}
The minimal generalized upper principal part $X_{G}^{U}$ is \textit{non-degenerate} if it satisfies all of the following conditions.
\begin{itemize}
    \item[(1)] The vector field $X$ is not free of characteristic orbits near infinity.
    \item[(2)] The polynomial vector field $X_{G}^{U}$ does not have curve of singularities intersecting the interior of $\mathbb{D}_{(\alpha,\beta)}$ and its boundary $\mathcal{D}$.
    \item[(3)] The infinity is \textit{non-dicritical} in the following sense. In the resolution of the singularities at infinity, the exceptional divisor does not contain arcs of singularities. 
\end{itemize}
\end{definition}

Let us give more details on the conditions of Definition \ref{def-mgupp-nd}. For this purpose, consider a polynomial vector field $X(x,y) = P(x,y)\frac{\partial}{\partial x} + Q(x,y)\frac{\partial}{\partial y}$ of homogeneous degree $\delta - 1$. We perform \textit{Bendixson compactification} (see \cite[Chapter 5.4]{DLA06}) by means of the coordinate change \begin{equation*}
\tilde{x}=\frac{x}{x^2+y^2},
\qquad
\tilde{y}=\frac{y}{x^2+y^2}.   
\end{equation*}

\begin{figure}[ht]\center{
\begin{overpic}[width=0.35\textwidth]{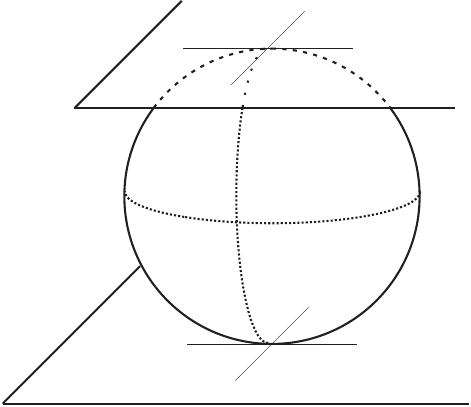}
%\begin{overpic}[grid,tics=10,width=0.35\textwidth]{fig-bendixson-compact.pdf}
\put(75,75){\footnotesize{$x$}}
\put(65,85){\footnotesize{$y$}}

\put(76,12){\footnotesize{$\tilde{x}$}}
\put(66,22){\footnotesize{$\tilde{y}$}}

\end{overpic}}
\caption{\footnotesize{Bendixson compactification of a vector field. the infinity of the phase space of $X$ (in $(x,y)$ coordinates) is the origin of the phase space of $\mathcal{B}(X)$ (in $(\tilde{x},\tilde{y})$ coordinates).}}
\label{fig-bendixson-compact}
\end{figure}
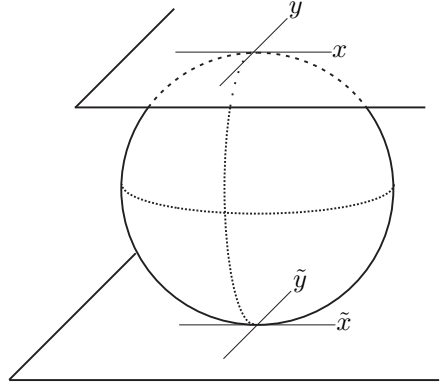

After multiplication by $({\tilde{x}^2+\tilde{y}^2})^{\delta-1}$ one obtains the vector field
\begin{equation}
 \begin{array}{rcl}
\mathcal{B}(X)(\tilde{x},\tilde{y})  & = & (\tilde{x}^2+\tilde{y}^2)^{\delta - 1}
\left[
(\tilde{y}^2-\tilde{x}^2)
P\!\left(
\frac{\tilde{x}}{\tilde{x}^2+\tilde{y}^2},
\frac{\tilde{y}}{\tilde{x}^2+\tilde{y}^2}
\right)
-2\tilde{x}\tilde{y}
Q\!\left(
\frac{\tilde{x}}{\tilde{x}^2+\tilde{y}^2},
\frac{\tilde{y}}{\tilde{x}^2+\tilde{y}^2}
\right)
\right]\frac{\partial}{\partial \tilde{x}}, \\
  & + & (\tilde{x}^2+\tilde{y}^2)^{\delta - 1}
\left[
(\tilde{x}^2-\tilde{y}^2)
Q\!\left(
\frac{\tilde{x}}{\tilde{x}^2+\tilde{y}^2},
\frac{\tilde{y}}{\tilde{x}^2+\tilde{y}^2}
\right)
-2\tilde{x}\tilde{y}
P\!\left(
\frac{\tilde{x}}{\tilde{x}^2+\tilde{y}^2},
\frac{\tilde{y}}{\tilde{x}^2+\tilde{y}^2}
\right)
\right]\frac{\partial}{\partial \tilde{y}},
\end{array}      
\end{equation}
where $\mathcal{B}(X)$ is the Bendixson compactification of $X$. Thus, the infinity of the phase space of $X$ is the origin of the phase space of $\mathcal{B}(X)$. See Figure \ref{fig-bendixson-compact}.

Now, the three conditions of Definition \ref{def-mgupp-nd} are equivalent to the following:
\begin{itemize}
    \item[(1)] The origin of the phase space of $\mathcal{B}(X)$ has characteristic orbit, that is, it is not a monodromic singularity.
    \item[(2)] The origin of the phase space of $\mathcal{B}(X_{G}^{U})$ is an isolated singularity, that is, there is no curve of singularities passing through the origin.
    \item[(3)] The origin of the phase space of $\mathcal{B}(X_{G}^{U})$ is not a dicritical singularity. 
\end{itemize}

The reader may notice that conditions (1) and (2) of Definition \ref{def-mgupp-nd} were also required in \cite[Theorem A]{DOP24}. On the other hand, in \cite{DOP24} the third condition required in order to assure topological equivalence near infinity was the Newton non-degeneracy. In the present paper, we allow the vector field to be Newton degenerate. However, when dropping the Newton non-degeneracy, with our approach we must require that the infinity is non dicritical in the sense of Definition \ref{def-mgupp-nd} item (3). The reason for this is that our approach relies in the Normal Form Theorem stated in Appendix \ref{appendix-nft} and the references therein, which concerns isolated singularities.

\begin{remark}\label{rem-dicritical}
It is possible to verify Definition \ref{def-mgupp-nd} item (3) in terms of Newton polytope. Indeed, it is equivalent to require that the polynomial in the left-hand side of Equation \eqref{eq-root-general-case} is not identically zero. However, in order to define such polynomial, it is necessary to introduce a lot of tools, which would postpone the statement of Theorem \ref{mthm} a lot. The reason why such polynomial appears will be clear in the proof of Theorem \ref{mthm} in Section \ref{sec-proofs}. Finally, we remark that, in the particular case where the upper diagram $\mathcal{P}^{U}$ has only one useful segment, Definition \ref{def-mgupp-nd} item (3) is equivalent to require that the polynomial in the left-hand side of Equation \eqref{eq-proof-roots} is not identically zero, which is easier to verify in comparison with the general case. 
\end{remark}

In what follows, we precisely define what we mean by topological equivalence near infinity.

\begin{definition}[Definition 10, \cite{DOP24}]\label{def-eq-infinity}
The planar polynomial vector fields $Y_{1}$ and $Y_{2}$ are \textit{topologically equivalent near infinity} if there is a homeomorphism $H:W_{1}\rightarrow W_{2}$ satisfying
\begin{itemize}
    \item[\textbf{(H1)}] $W_{1}$ and $W_{2}$ are open sets containing the boundary $\mathcal{D}$ of the PL-disk $\mathbb{D}_{(\alpha,\beta)}$, and $H(\mathcal{D}) = \mathcal{D}$.
    \item[\textbf{(H2)}] Denote the PL compactification of $Y_{i}$ by $\overline{Y}_{i}$, for $i = 1,2$.  Given $t > 0$ and $p\in W_{1}$, there is $t^{*} > 0$ such that
\begin{equation*}
H\left(\varphi_{\overline{Y}_{1}}\left(p,[0,t]\right)\right) = \varphi_{\overline{Y}_{2}}\left(H\left(p\right),[0,t^{*}]\right),    
\end{equation*}
in which $\varphi_{\overline{Y}_{i}}$ is the flow of $\overline{Y}_{i}$, for $i = 1,2$.
\end{itemize}
\end{definition}

Now, we are in position to state the main result of this paper.

\begin{mtheorem}\label{mthm}
Fix a Newton polytope $\mathcal{P}$ and consider the set $\mathfrak{U}_{1}(\mathcal{P})$. There is a subset $\widetilde{\mathfrak{U}}_{1}(\mathcal{P})\subset \mathfrak{U}_{1}(\mathcal{P})$, which is open and dense in $\mathfrak{U}_{1}(\mathcal{P})$, satisfying the following property: If $X\in \widetilde{\mathfrak{U}}_{1}(\mathcal{P})$ and its minimal generalized upper principal part $X_{G}^{U}$ is non-degenerate, then $X$ and $X_{G}^{U}$ are topologically equivalent near infinity. 
\end{mtheorem}

In Theorem \ref{mthm}, the set $\widetilde{\mathfrak{U}}_{1}(\mathcal{P})$ is a Zariski-open subset of $\mathfrak{U}_{1}(\mathcal{P})$. The proof of Theorem \ref{mthm} is given in Section \ref{sec-proofs} using the so called \textit{compactification adapted to $\mathcal{P}$} given in \cite[Section 3]{DOP24}, and a proof for the case which the upper diagram has only one useful segment is given in Section \ref{sec-one-useful-segment} using Poincaré--Lyapunov compactification. For the sake of readability, firstly we prove the particular case and afterwards we prove the general case. Indeed, once the proof for such particular case is well understood, it is easier to follow the proof of the general case. 

\begin{remark}
The set $\widetilde{\mathfrak{U}}_{1}(\mathcal{P})$ is algebraically defined in Equations \eqref{eq-def-u1} and \eqref{eq-def-u1-final} in the proof of Theorem \ref{mthm} in Section \ref{sec-proofs}. However, analogously as in Remark \ref{rem-dicritical}, in order to precisely define Equations \eqref{eq-def-u1} and \eqref{eq-def-u1-final} it is necessary to introduce a lot of tools that would postpone the statement of Theorem \ref{mthm} a lot. Finally, it is important remark is that the subsets $\widetilde{\mathfrak{U}}_{1}(\mathcal{P})$ and $\mathfrak{U}_{2}(\mathcal{P})$ of $\mathfrak{U}_{1}(\mathcal{P})$ do not depend on the initial vector field $X$.
\end{remark}

The following Lemma, which is a particular case of \cite[Proposition 4]{Zup96} (see also \cite{Zup00}), play a crucial role in the proof or our main result.

\begin{lemma}\label{prop-normal-form-semihyp}
Consider the planar analytic vector field
\begin{equation}\label{eq-vf-prop}
X(x,y) = \left(\sum_{r + s \geq 1}a_{r,s}x^{r}y^{s}\right)\frac{\partial}{\partial x} + \left(\sum_{r + s \geq 1}b_{r,s}x^{r}y^{s}\right)\frac{\partial}{\partial y},
\end{equation}
and suppose that the origin is a semi-hyperbolic singularity, and the linearization of $X$ satisfies $a_{1,0} \neq 0$. Denote by $l$ the smallest positive integer such that $X$ has the monomial $a_{0,l}y^{l}\frac{\partial}{\partial x}$ or $b_{0,l+1}y^{l+1}\frac{\partial}{\partial y}$. Then there exists a function $F$, which is polynomial with respect to $(b_{0,l+1}, a_{0,l}, b_{1,1})$, such that if the coefficients $(b_{0,l+1}, a_{0,l}, b_{1,1})$ of the vector field $X$ given in \eqref{eq-vf-prop} satisfy $F(b_{0,l+1}, a_{0,l}, b_{1,1}) \neq 0$, then $X$ is topologically equivalent near the origin to
\begin{equation*}
Y(x,y) = a_{1,0}x \frac{\partial}{\partial x} + F\left(b_{0,l+1}, a_{0,l}, b_{1,1}\right)y^{l+1} \frac{\partial}{\partial y}. 
\end{equation*}

Moreover, the polynomial $F$ is given by:
\begin{itemize}
    \item[(1)] If $l = 1$ and $b_{2,0} = 0$ in Equation \eqref{eq-vf-prop}, then $F(b_{0,2}, a_{0,1}, b_{1,1}) = b_{0,2} - \frac{a_{0,1}b_{1,1}}{a_{1,0}}$.
    \item[(2)] If $l = 2$, then $F(b_{0,3}, a_{0,2}, b_{1,1}) = b_{0,3}$.
    \item[(3)] If $l \geq 3$, then $F(b_{0,l+1}, a_{0,l}, b_{1,1}) = b_{0,l+1} - \frac{a_{0,l}b_{1,1}}{a_{1,0}}$.
\end{itemize}
\end{lemma}

The proof of Lemma \ref{prop-normal-form-semihyp} is technical, and it is a consequence of the Normal Form Theorem stated in Appendix \ref{appendix-nft}. For the sake of readability, we postpone the proof of Lemma \ref{prop-normal-form-semihyp} to Section \ref{proof-prop}. It is important to keep in mind that the topological normal form will depend on $a_{1,0}x\frac{\partial}{\partial x}$, $b_{0,l+1}y^{l+1}\frac{\partial}{\partial y}$, $a_{0,l}y^{l}\frac{\partial}{\partial x}$ and $b_{1,1}xy\frac{\partial}{\partial y}$, with $l$ being the smallest positive integer such that $X$ has the monomial $a_{0,l}y^{l}\frac{\partial}{\partial x}$ or $b_{0,l+1}y^{l+1}\frac{\partial}{\partial y}$.

In Item (1) of Lemma \ref{prop-normal-form-semihyp}, if one would not require the condition $b_{2,0} = 0$, then the expression of the polynomial $F$ would be slightly different. However, for our purposes, it is enough to consider the case where $b_{2,0} = 0$. The reason why will be clear in the proof of Theorem \ref{mthm} in Section \ref{sec-proofs} (and also in Section \ref{sec-one-useful-segment}, for the particular case).

\section{The case with one useful segment}\label{sec-one-useful-segment}

In this section, we assume that the upper principal part has only one useful segment with negative slope, which will be simply denoted by $\gamma^{U}$. Observe that with this assumption the upper diagram $\mathcal{P}^{U}$ can have at most three compact segments. The purpose of this section is to provide a proof of Theorem \ref{mthm} in the case where the upper principal part has only one useful segment. In the proof of this particular case, we use Poincaré--Lyapunov compactification. 

\begin{figure}[ht]\center{
\begin{overpic}[width=0.8\textwidth]{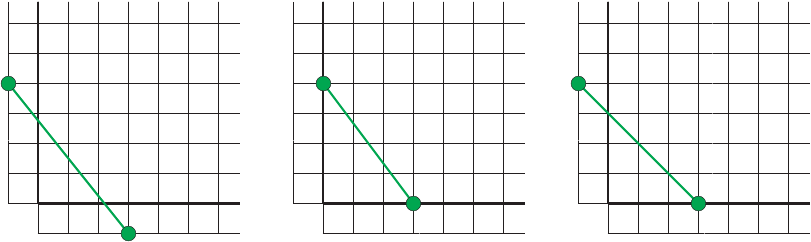}
%\begin{overpic}[grid,tics=10,width=0.7\textwidth]{fig-cases-polygon.pdf}
\put(-15,20){\footnotesize{$(- 1, \alpha -1)$}}
\put(10,-2){\footnotesize{$(\beta- 1,  -1)$}}

\put(30,20){\footnotesize{$(0, \alpha )$}}
\put(50,-2){\footnotesize{$(\beta,  0)$}}

\put(62,20){\footnotesize{$(-1, \alpha )$}}
\put(82,-2){\footnotesize{$(\beta -1,  0)$}}

\put(10,31){\footnotesize{Case I}}
\put(50,31){\footnotesize{Case II}}
\put(82,31){\footnotesize{Case III}}
\end{overpic}}
\caption{\footnotesize{Possible cases of the unique useful segment $\gamma^{U}$, which is highlighted in green.}}
\label{fig-cases-polygon}
\end{figure}
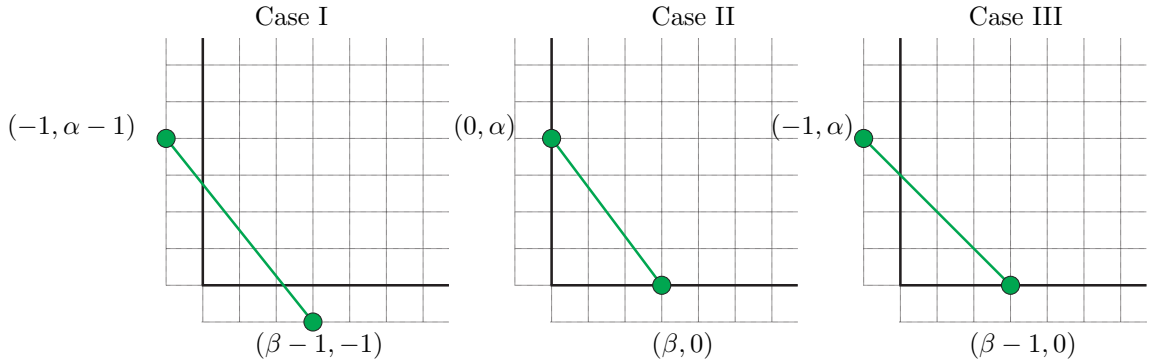

Once the useful segment $\gamma^{U}$ is fixed, we choose the weight vector $\omega$ for the Poincaré--Lyapunov compactification as the normal vector of the support line of $\gamma^{U}$, and we will denote it by $\omega = (\alpha, \beta)$. Here, $\operatorname{gcd}(\alpha,\beta) = 1$ and $\alpha,\beta > 0$. This implies that we have three possibilities, which are sketched in Figure \ref{fig-cases-polygon}.

\hspace{1cm}

For each case in Figure \ref{fig-cases-polygon}, the support line of the useful side $\gamma^{U}$ is given by
\begin{itemize}
    \item $\{(m,n); \ \alpha m + \beta n = \delta\}$ with $\delta = \alpha\beta - \alpha - \beta$, in Case I. Here, we do not have non useful segments.
    \item $\{(m,n); \ \alpha m + \beta n = \delta\}$ with $\delta = \alpha\beta$, in Case II. Here, one may have at most two non useful segments.
    \item $\{(m,n); \ \alpha m + \beta n = \delta\}$ with $\delta = \alpha\beta - \alpha$, in Case III. Here, one may have at most one non useful segment contained in $\mathbb{Z}\times [-1,0]$. 
\end{itemize}

\begin{remark}
One would have a fourth case, in which the useful segment $\gamma^{U}$ contains the points $(0,\alpha -1)$ and $(\beta, -1)$. However, this case is equivalent to Case III. Indeed, it is just a matter of applying a change of coordinates of the form $(x,y) \mapsto (\tilde{y},\tilde{x})$.    
\end{remark}

Consider the support line of the set $\Gamma^{(1)}_{1}$ (recall the definition of such set in Section \ref{sec-def-gupp}). Such a support line below $\gamma^{U}$ has the form $\tilde{r} = \{(m,n); \ \alpha m + \beta n = \widetilde{\delta}\}$ with $-1 \leq \widetilde{\delta} < \delta$. Therefore, the minimal generalized upper principal part is $X_{G}^{U} = X_{\Delta}^{U} + X_{1}^{U}$, in which the monomials of $X_{1}^{U}$ are related to points of the support $\mathcal{S}$ that are contained in the set $\Gamma^{(1)}$. Recall that the minimal generalized upper principal part $X_{G}^{U}$ is non-degenerate if it satisfies Definition \ref{def-mgupp-nd}.

By hypothesis, $\gamma^{U}$ is the only useful segment and it has negative slope. Therefore, after performing Poincaré--Lyapunov compactification, the origin in both positive and negative $y$-directions will be either an elementary singularity or a regular point. In the same fashion, since the upper diagram has only one useful segment, after performing Poincaré--Lyapunov compactification, the origin in both positive and negative $x$-directions will be either an elementary singularity or a regular point. These facts follow from \cite[Propositions 13 and 14]{DOP24}). Therefore, in each chart of the compactification, it is sufficient to study the singularities along the curve at infinity but outside the origin. Computations in both positive and negative $x$-directions are analogous, so in what follows we perform the computations in the positive $x$-direction.

\begin{remark}
As proved in \cite[Propositions 15]{DOP24}, the Newton non-degeneracy implies that, at infinity, all singularities outside the origin of the directional charts are elementary (that is, they are hyperbolic or semi-hyperbolic). However, in the present paper, we drop the Newton non-degeneracy condition, and therefore it might be possible that non-elementary singularities at infinity appear outside the origin.
\end{remark}

\subsection{The expression of the compactified vector field}

Firstly, let us introduce some notation. The vector field $X_{\Gamma_{1}^{(0)}}^{U} = (P_{0}^{U},Q_{0}^{U})$ has monomials related to points in the set $\Gamma_{1}^{(0)}\subset \gamma^{U}$ and $X_{\Gamma_{1}^{(1)}}^{U} = (P_{1}^{U},Q_{1}^{U})$ has monomials related to points in the set $\Gamma_{1}^{(1)}$. In general, $X_{\Delta}^{U}$ has more monomials than $X_{\Gamma_{1}^{(0)}}^{U}$. %and the vector field $X_{G}^{U}$ is the sum of $X_{\Gamma_{1}^{(0)}}^{U} + X_{\Gamma_{1}^{(1)}}^{U}$ plus one monomial whose point is in $\{-1\}\times \mathbb{Z}$ and one monomial whose point is in $\mathbb{Z} \times \{-1\}$ (when such points of the support exist).

Since both $X_{\Gamma_{1}^{(0)}}^{U} = (P_{0}^{U},Q_{0}^{U})$ and $X_{\Gamma_{1}^{(1)}}^{U} = (P_{1}^{U},Q_{1}^{U})$ are quasi-homogeneous vector fields of type $\omega = (\alpha, \beta)$ and degree $\delta$ and $\widetilde{\delta}$, respectively, they satisfy \begin{multicols}{2}\noindent
\begin{equation*}
P_{0}^{U}(c^{\alpha}x,c^{\beta}y) = c^{\delta + \alpha}P_{0}^{U}(x,y),
\end{equation*}
\begin{equation*}
P_{1}^{U}(c^{\alpha}x,c^{\beta}y) = c^{\widetilde{\delta} + \alpha}P_{1}^{U}(x,y),
\end{equation*}
\begin{equation*}
Q_{0}^{U}(c^{\alpha}x,c^{\beta}y) = c^{\delta + \beta}Q_{0}^{U}(x,y),
\end{equation*}
\begin{equation*}
Q_{1}^{U}(c^{\alpha}x,c^{\beta}y) = c^{\widetilde{\delta} + \beta}Q_{1}^{U}(x,y),
\end{equation*}
\end{multicols}
\noindent
and then the compactification of $X$ in the positive $x$-direction will have the form
\begin{equation}\label{eq-proof-compactified-x}
\left\{
  \begin{array}{rl}
   \dot{u} & = \left( Q_{0}^{U}(1,u) - \displaystyle\frac{\beta}{\alpha}uP_{0}^{U}(1,u)\right) + v^{\delta - \widetilde{\delta}}\left(Q_{1}^{U}(1,u) - \displaystyle\frac{\beta}{\alpha}uP_{1}^{U}(1,u) \right) + O(v^{\delta - \widetilde{\delta} + 1}), \\
   \dot{v} & = -\displaystyle\frac{v}{\alpha}\left( P_{0}^{U}(1,u) + v^{\delta - \widetilde{\delta}}P_{1}^{U}(1,u) + O(v^{\delta - \widetilde{\delta} + 1})\right).
  \end{array}
\right.    
\end{equation}

Observe that higher order terms in Equation \eqref{eq-proof-compactified-x} are related to points in $\Gamma_{1}^{(i)}$, with $i\geq 2$. In addition, singularities of \eqref{eq-proof-compactified-x} at infinity $\{v = 0\}$ (if they exist) are points of the form $(\lambda,0)$, in which $\lambda$ is a real solution of the polynomial equation
\begin{equation}\label{eq-proof-roots}
Q_{0}^{U}(1,u) - \displaystyle\frac{\beta}{\alpha}uP_{0}^{U}(1,u) = 0. 
\end{equation}

Due to item (3) of Definition \ref{def-mgupp-nd}, we must suppose that the polynomial $Q_{0}^{U}(1,u) - \frac{\beta}{\alpha}uP_{0}^{U}(1,u)$ is not identically zero, otherwise one would have curves of singularities contained at infinity (see also Remark \ref{rem-dicritical}). Since the upper diagram has only one useful segment, then the origin $\lambda = 0$ is hyperbolic or semi-hyperbolic (if it is a singularity). Indeed, the origin is either a hyperbolic singularity, or it is semi-hyperbolic whose central manifolds are tangent to $\{v = 0\}$, and the (un)stable manifold is transversal to $\{v = 0\}$. In any case, the topological behavior of the singularity depends on monomials related to $\gamma^{U}$.

The next step is to study the singularities at infinity but outside the origin of the positive $x$ and $y$ directional charts.

\begin{lemma}\label{lemma-compact-x-trans}
If $\lambda\in\mathbb{R}$ is a non zero root of \eqref{eq-proof-roots}, we perform a translation of the form 
$$u = z + \lambda, \quad v = v,$$
and then one obtains the system
\begin{equation}\label{eq-compact-x-trans}
\left\{
  \begin{array}{rcrcrcc}
   \dot{z} & = & \displaystyle\sum_{k = 1}^{\alpha + 1}A_{k}^{(0)}z^{k} & + &  v^{\delta - \widetilde{\delta}}\displaystyle\sum_{k = 0}^{N_{1}}A_{k}^{(1)}z^{k} & + & O(v^{\delta - \widetilde{\delta} + 1}), \\
   \dot{v} & = & -v\displaystyle\sum_{k = 0}^{\alpha + 1}B_{k}^{(0)}z^{k} & -&  v^{\delta - \widetilde{\delta}+1}\displaystyle\sum_{k = 0}^{N_{2}}B_{k}^{(1)}z^{k} & + & O(v^{\delta - \widetilde{\delta} + 2}),
  \end{array}
\right.   
\end{equation} 
whose coefficients are given by
\begin{equation*}
\begin{array}{ccccc}
   A_{k}^{(0)} & = & A_{k}^{(0)}(\lambda) & = & \displaystyle\frac{1}{k!} \left( \displaystyle\frac{\partial^{k}Q_{0}^{U}}{\partial y^{k}}(1,\lambda) - \frac{\beta}{\alpha} \left( k\displaystyle\frac{\partial^{k-1}P_{0}^{U}}{\partial y^{k-1}}(1,\lambda) + \lambda\displaystyle\frac{\partial^{k}P_{0}^{U}}{\partial y^{k}}(1,\lambda) \right)\right),  \\
   & & \\
   A_{k}^{(1)} & = & A_{k}^{(1)}(\lambda) & = & \displaystyle\frac{1}{k!} \left( \displaystyle\frac{\partial^{k}Q_{1}^{U}}{\partial y^{k}}(1,\lambda) - \frac{\beta}{\alpha} \left( k\displaystyle\frac{\partial^{k-1}P_{1}^{U}}{\partial y^{k-1}}(1,\lambda) + \lambda\displaystyle\frac{\partial^{k}P_{1}^{U}}{\partial y^{k}}(1,\lambda) \right)\right),
\end{array}
\end{equation*}
\begin{equation*}
\begin{array}{ccccccccccc}
B_{k}^{(0)} & = & B_{k}^{(0)}(\lambda) & = & \displaystyle\frac{1}{\alpha}\frac{1}{k!}\frac{\partial^{k}P_{0}^{U}}{\partial y^{k}}(1,\lambda), & \ & B_{k}^{(1)} & = & B_{k}^{(1)}(\lambda) & = &  \displaystyle\frac{1}{\alpha}\frac{1}{k!}\frac{\partial^{k}P_{1}^{U}}{\partial y^{k}}(1,\lambda).
\end{array}
\end{equation*}
\end{lemma}

For the sake of readability, we leave the proof of Lemma \ref{lemma-compact-x-trans} to Appendix \ref{appendix-proofs}.

\begin{remark}
In Lemma \ref{lemma-compact-x-trans}, the coefficients $A_{k}^{(0)}$ and $B_{k}^{(0)}$ depend on the monomials related to points of $\Gamma_{1}^{(0)}\subset\gamma^{U}$, whereas $A_{k}^{(1)}$ and $B_{k}^{(1)}$ depend on monomials related to the set $\Gamma_{1}^{(1)}$ below $\gamma^{U}$. On the other hand, the index $k$ is related to powers of $z$ in Equation \eqref{eq-compact-x-trans}. In particular, the coefficients $A_{1}^{(0)}$, $B_{0}^{(0)}$, $B_{1}^{(0)}$, $A_{0}^{(1)}$ and $B_{0}^{(1)}$ are given by
$$
A_{1}^{(0)}(\lambda) = \displaystyle\frac{\partial Q_{0}^{U}}{\partial y}(1,\lambda) - \frac{\beta}{\alpha}\left( P_{0}^{U}(1,\lambda) + \lambda\displaystyle\frac{\partial P_{0}^{U}}{\partial y}(1,\lambda) \right), \quad A_{0}^{(1)}(\lambda) = Q_{1}^{U}(1,\lambda) - \displaystyle\frac{\beta}{\alpha}\lambda P_{1}^{U}(1,\lambda),
$$
$$
B_{0}^{(0)}(\lambda) = \frac{1}{\alpha}P_{0}^{U}(1,\lambda), \quad B_{1}^{(0)}(\lambda) = \frac{1}{\alpha}\frac{\partial P_{0}^{U}}{\partial y}(1,\lambda),  \quad B_{0}^{(1)}(\lambda) = \frac{1}{\alpha}P_{1}^{U}(1,\lambda).
$$
\end{remark}

\subsection{Computing the topological normal form}\label{sec-nf-one-seg}
The next step is to compute the $C^{0}$ normal form of System \eqref{eq-compact-x-trans}. One should compute the normal form in order to define the open and dense subset in $\widetilde{\mathfrak{U}}_{1}(\mathcal{P})\subset\mathfrak{U}_{1}(\mathcal{P})$. Observe that higher order terms in Equation \eqref{eq-compact-x-trans} are related to points in $\Gamma_{1}^{(i)}$, with $i\geq 2$. In what follows, we pay attention to the monomials
\begin{equation*}
A_{1}^{(0)}z\frac{\partial}{\partial z}, \quad A_{0}^{(1)}v^{\delta - \widetilde{\delta}}\frac{\partial}{\partial z}, \quad B_{0}^{(0)}v\frac{\partial}{\partial v}, \quad B_{1}^{(0)}zv\frac{\partial}{\partial v}, \quad B_{0}^{(1)}v^{\delta - \widetilde{\delta} + 1}\frac{\partial}{\partial v}.
\end{equation*} 

We are interested in the cases in which the singularities along the infinity are all elementary, that is, they are either hyperbolic or semi-hyperbolic. Observe that, if the coefficient $B_{0}^{(0)}(\lambda)$ in Lemma \ref{lemma-compact-x-trans} is different from zero, then the singularity $\lambda$ is either hyperbolic or semi-hyperbolic. In the former case, the topological normal form of this singularity is given by the Grobman--Hartman Theorem. Otherwise, it would be a semi-hyperbolic singularity whose center manifolds are tangent to the line at infinity $\{v = 0\}$, and the (un)stable manifold is transversal to infinity. In both cases, the behavior of the vector field near the singularity $\lambda$ (under topological equivalence) depends only on the monomials associated to the set $\Gamma_{1}^{(0)}\subset\gamma^{U}$. In summary, if $(\lambda,0)$ is a singularity at infinity with $B_{0}^{(0)}(\lambda) \neq 0$, then the Poincaré--Lyapunov compactifications of $X$ and $X_{G}^{U}$ are topologically equivalent near $(\lambda,0)$.

Now, one must consider the case $B_{0}^{(0)}(\lambda) = 0$ and $A_{1}^{(0)}(\lambda) \neq 0$, so $\lambda$ is a semi-hyperbolic singularity. In this case, the (un)stable manifold is tangent to the line at infinity $\{v = 0\}$, whereas center manifolds are not. Therefore, one should add extra conditions in the coefficients of \eqref{eq-compact-x-trans} in order to assure the topological equivalence of $X$ and $X^{U}_{G}$. 

\begin{remark}
Observe that, in the case $B_{0}^{(0)}(\lambda) = A_{1}^{(0)}(\lambda) = 0$, the singularity $(\lambda,0)\in\{v = 0\}$ is neither hyperbolic nor semi-hyperbolic, therefore it is non-elementary. This would lead to further blow-up analysis and we do not treat this case in this paper.  
\end{remark}

Observe that System \eqref{eq-compact-x-trans} fits the assumptions of Lemma \ref{prop-normal-form-semihyp}. Then, it follows that there is a polynomial function $F_{\lambda}$ such that System \eqref{eq-compact-x-trans} is topologically equivalent to system
\begin{equation*}\label{eq-compact-x-nf}
\left\{
  \begin{array}{rcl}
   \dot{z} & = & A_{1}^{(0)}z, \\
   \dot{v} & = & F_{\lambda}\left(B_{0}^{(1)},A_{0}^{(1)},B_{1}^{(0)}\right)v^{\delta - \widetilde{\delta}+1},
  \end{array}
\right.   
\end{equation*}
provided that $F_{\lambda}\left(B_{0}^{(1)},A_{0}^{(1)},B_{1}^{(0)}\right)\neq 0$.

For each singularity $(\lambda,0)$ at infinity, define the set
\begin{equation*}
 \mathfrak{U}_{2,\lambda}(\mathcal{P}) := \left\{ X\in\mathfrak{U}_{1}(\mathcal{P}); \quad B_{0}^{(0)} = 0 \right\} \cap \left\{X\in\mathfrak{U}_{1}(\mathcal{P}); \quad F_{\lambda}\left(B_{0}^{(1)},A_{0}^{(1)},B_{1}^{(0)}\right) = 0 \right\}.   
\end{equation*}

It follows from our reasoning that if $X\not\in \mathfrak{U}_{2,\lambda}(\mathcal{P})$, then the compactifications of $X$ and $X_{G}^{U}$ are topologically equivalent near $(\lambda,0)$. We further define
\begin{equation*}
 \mathfrak{U}_{2}(\mathcal{P}) := \bigcup_{\lambda\in\mathcal{Z}} \mathfrak{U}_{2,\lambda}(\mathcal{P}),  \qquad  \widetilde{\mathfrak{U}}_{1}(\mathcal{P}) := \mathfrak{U}_{1}(\mathcal{P})\backslash \mathfrak{U}_{2}(\mathcal{P}),
\end{equation*}
where $\mathcal{Z}$ is the set of all real roots of the  polynomial Equation \eqref{eq-proof-roots} (that is, $\lambda$ is a singularity at infinity $\{v = 0\}$). Of course, $\mathcal{Z}$ is a finite set. Thus, if $X\in\widetilde{\mathfrak{U}}_{1}(\mathcal{P})$, then the compactifications of $X$ and $X_{G}^{U}$ are topologically equivalent near each singularity $(\lambda,0)$ at infinity. The set $\widetilde{\mathfrak{U}}_{1}(\mathcal{P})$ is open and dense in $\mathfrak{U}_{1}(\mathcal{P})$ with respect to the Zariski topology. Finally, observe that the subsets $\widetilde{\mathfrak{U}}_{1}(\mathcal{P})$ and $\mathfrak{U}_{2}(\mathcal{P})$ of $\mathfrak{U}_{1}(\mathcal{P})$ do not depend on the initial vector field $X$ given, because each set $\mathfrak{U}_{2,\lambda}(\mathcal{P})$ is defined independently of $X$.

\subsection{The existence of the topological equivalence}\label{sec-existence-top-eq}

The proof so far guarantees that the compactifications of $X$ and $X_{G}^{U}$ are topologically equivalent in a neighborhood of each singularity $(\lambda,0)$ at infinity. The next step is to prove that the compactifications of both $X$ and $X_{G}^{U}$ are topologically equivalent in a neighborhood of the whole boundary of the PL-disk. This can be done in the very same fashion as in \cite[Theorem B]{Dum77} and \cite[Theorem A]{DOP24} (see \cite[Theorem A]{PerSil22} for an analogous construction in the context of topological classification of singularities of constrained differential systems).

We briefly sketch the proof that assures the existence of a topological equivalence between the compactifications of $X$ and $X_{G}^{U}$ in a neighborhood of the boundary of the PL-disk. Due to our non degeneracy assumptions, there is a finite number of singularities at infinity, all of them hyperbolic or semi-hyperbolic, and the infinity is homeomorphic to $\mathbb{S}^{1}$. The infinity is decomposed in a finite number of sectors, and in each sector (the compactifications of) $X$ and $X_{G}^{U}$ are topologically equivalent. Applying the Pasting Lemma, one can “glue” two homeomorphisms defined in two adjacent sectors. Repeating this reasoning to all sectors, we obtain a homeomorphism that gives the topological equivalence defined in a neighborhood of the whole infinity. It is important to remark that, since the infinity is homeomorphic to $\mathbb{S}^{1}$, in the end of the construction one must glue the homeomorphisms defined in the last and the first sector. This can be done under the hypothesis that at least one singularity at infinity has characteristic orbit, that is, one must avoid a monodromic behavior near infinity.

%Moreover, the term $\widetilde{b}_{0,m}x^{m}$ satisfies 
%\begin{equation*}
%m = \operatorname{min}\left\{\begin{array}{ccl}
%   k,  & \text{where} & b_{0,k} \neq 0, \\
%   lr + s  & \text{where} & b_{r,s} \neq 0 \text{ and } a_{0,l} \neq 0,
%\end{array}\right.  \quad
%\widetilde{b}_{0,m} = \left\{\begin{array}{ccl}
%   b_{0,k},  & \text{if} & m = k, \\
%   f(b_{r,s}, a_{0,l})  & \text{if} & m = lr + s, \\
%   g(b_{0,k}, b_{r,s}, a_{0,l})& \text{if} & m = k = lr + s,
%\end{array}\right.
%\end{equation*}
%for some polynomial functions $f$ and $g$.
%\end{proposition}

%One consequence of Lemma \ref{prop-normal-form-semihyp} is the following.

%\begin{proposition}\label{prop-normal-form-semihyp-c0}
%The vector field \eqref{eq-nf-prop} is topologically equivalent ($C^{0}$-equivalent) to
%\begin{equation}\label{eq-nf-prop-c0}
%Y(x,y) = a_{1,0}x \frac{\partial}{\partial x} + \widetilde{b}_{0,m}y^{m} \frac{\partial}{\partial y},
%\end{equation}
%in which the exponent $m$ and the coefficients $a_{1,0}$ and $\widetilde{b}_{0,m}$ are given as in Lemma \ref{prop-normal-form-semihyp}.
%\end{proposition}

\section{Proof of the main Theorem in the general case}\label{sec-proofs}

This proof of Theorem \ref{mthm} for the general case is carried using the so called \textit{compactification adapted to $\mathcal{P}$} as described in \cite[Section 3]{DOP24}. Such compactification is essentially a toric compactification as described in \cite{Kho77, Kho78}. In what follows we recall the necessary definitions and we refer to \cite[Section 3]{DOP24} and the references therein for details on this compactification procedure.

\begin{definition}[Definition 2, \cite{DOP24}]
 Let $\Sigma^{*}_{\mathcal{P}} = \{(\alpha_{i},\beta_{i})\}_{i=0}^{l+1}$ be the collection of all inward normal vectors of $\mathcal{P}^{U}$ with $\gcd(\alpha_{i},\beta_{i}) = 1$. A \emph{simple fan} is the collection of vectors $\Sigma_{\mathcal{P}} = \{\xi_{j}\}_{j = 0}^{s}$, with $\xi_{j} = (\mu_{j},\nu_{j})\in\mathbb{Z}^{2}$, satisfying the following conditions.
\begin{itemize}
    \item[\textbf{(i)}] $\Sigma^{*}_{\mathcal{P}}\subset \Sigma_{\mathcal{P}}$;
    \item[\textbf{(ii)}] $\xi_{0} = (0,1)$ and $\xi_{s} = (1,0)$;
    \item[\textbf{(iii)}] $\gcd(\mu_{j},\nu_{j}) = 1$ and two consecutive vectors satisfy $\det\left(
  \begin{array}{cc}
    \mu_{j-1} & \nu_{j-1} \\
    \mu_{j} & \nu_{j} \\
  \end{array}
\right) = 1$;
    \item[\textbf{(iv)}] Two consecutive vectors $\xi_{j-1},\xi_{j}$ are not normal to two consecutive segments $\mathcal{P}^{U}$; 
    \item[\textbf{(v)}] $\Sigma_{\mathcal{P}}$ is \textit{minimal} in the sense that there is no collection satisfying conditions (i) - (iv) containing a smaller number of vectors.
\end{itemize}   
\end{definition}

Observe that there is no vector $\xi_{j}\in\Sigma_{\mathcal{P}}$ which is contained in the first quadrant of $\mathbb{R}^{2}$, for $j = 1,\dots, s-1$. For $j = 1,\dots,s$, each pair $B_{j} = \{\xi_{j-1},\xi_{j}\}$ defines a diffeomorphism $\psi_{j}$ as
\begin{equation*}
\begin{array}{cccccccc}
\psi_{j}: & (\mathbb{R}^{*})^{2} & \rightarrow &  (\mathbb{R}^{*})^{2} & ; & \psi_{j}(u,v) & = & (u^{\mu_{j-1}}v^{\mu_{j}}, \ u^{\nu_{j-1}}v^{\nu_{j}}); \\
\psi_{j}^{-1}: & (\mathbb{R}^{*})^{2} & \rightarrow &  (\mathbb{R}^{*})^{2} & ; & \psi_{j}^{-1}(x,y) & = & (x^{\nu_{j}}y^{-\mu_{j}}, \ x^{-\nu_{j-1}}y^{\mu_{j-1}}); \\
\end{array}
\end{equation*}
in which $\mathbb{R}^{*} = \mathbb{R}\backslash\{0\}$. In particular, the set $B_{0} = \{\xi_{s},\xi_{0}\} = \{(1,0), (0,1)\}$ induces the identity map. For $j = 1$ and $j = s$ we obtain the maps $\psi_{1}$ and $\psi_{s}$ as 
$$
\begin{array}{cccccccc}
\psi_{1}: & \mathbb{R}\times\mathbb{R}^{*} & \rightarrow &  \mathbb{R}^{*} \times \mathbb{R} & ; & \psi_{1}(u,v) & = & (v^{-1}, \ uv^{\nu_{1}}); \\
\psi_{1}^{-1}: & \mathbb{R}^{*}\times\mathbb{R} & \rightarrow &  \mathbb{R}\times\mathbb{R}^{*} & ; & \psi_{1}^{-1}(x,y) & = & (x^{\nu_{1}}y, \ x^{-1});
\end{array}
$$
$$
\begin{array}{cccccccc}
\psi_{s}: & \mathbb{R}^{*}\times\mathbb{R} & \rightarrow &  \mathbb{R}\times\mathbb{R}^{*} & ; & \psi_{s}(u,v) & = & (u^{\mu_{s-1}}v, \ u^{-1}); \\
\psi_{s}^{-1}: & \mathbb{R}\times\mathbb{R}^{*} & \rightarrow &  \mathbb{R}^{*} \times \mathbb{R} & ; & \psi_{s}^{-1}(x,y) & = & (y^{-1}, \ xy^{\mu_{s-1}}).
\end{array}
$$

\begin{definition}[Definition 3, \cite{DOP24}]
A \emph{compactification of $(\mathbb{R}^{*})^{2}$ adapted to $\mathcal{P}$} is a pair $(\mathcal{M},\Psi)$ satisfying the following conditions.
\begin{itemize}
    \item[\textbf{(a)}] $\mathcal{M}$ is a real smooth compact 2-dimensional manifold containing $(\mathbb{R}^{*})^{2}$;
    \item[\textbf{(b)}] $\Psi:\mathcal{M}\rightarrow\mathbb{R}^{2}$ is a map locally given by $\overline{\psi}_{j}:U_{j}\subset\mathcal{M}\rightarrow\mathbb{R}^{2}$, which is the extension of $\psi_{j}$ to a set $\mathcal{I}$ that is the union of 1-dimensional smooth manifolds in general position;
    \item[\textbf{(c)}] The collection  $U_{j}$ covers $\mathcal{M}$.
\end{itemize}
\end{definition}

\begin{definition}[Definition 4, \cite{DOP24}]\label{def-compact}
A \emph{compactification of the polynomial vector field $X$ adapted to $\mathcal{P}$} is the vector field $\overline{X}:\mathcal{M}\rightarrow T\mathcal{M}$ defined as the pullback of $X$ by $\Psi$, that is, $\overline{X} = \Psi^{*}X$. See Figure \ref{fig-diag-def-compact}.
\end{definition}

\begin{figure}[h!]
\begin{flushright}
\begin{center}
\begin{tikzpicture}
\node (A) {$T\mathcal{M}$};
\node (B) [right of=A] {$T\mathbb{R}^{2}$};
\node (C) [below of=A] {$\mathcal{M}$};
\node (D) [below of=B] {$\mathbb{R}^{2}$};
\large\draw[->] (C) to node {\mbox{{\footnotesize $\overline{X}$}}} (A);
%\draw[->] (X) to node [swap] {\mbox{{\footnotesize $q$}}} (Z);
\large\draw[->] (D) to node {\mbox{{\footnotesize $X$}}} (B);
\large\draw[->] (A) to node {\mbox{{\footnotesize $\Psi^{*}$}}} (B);
\large\draw[->] (C) to node {\mbox{{\footnotesize $\Psi$}}} (D);
\end{tikzpicture}
\end{center}
\end{flushright}
\caption{\footnotesize{Cummutative diagram of Definition \ref{def-compact}.}}
\label{fig-diag-def-compact}
\end{figure}

Consider the vector field $X_{j} = \psi_{j}^{*} X$, that is, $X_{j}$ is the pushforward of $X$ by $\psi_{j}$, for $j = 1,\dots,s$. Observe that $\psi_{j}$ is a diffeomorphism outside $\{uv = 0\}$ and $X_{j}$ can be analytically extended to such set. Such extension will be denoted by $\overline{X}_{j}$, and the domain of $\overline{X}_{j}$ will be denoted by $U_{j}$, where $\{uv = 0\}$ plays the role of infinity. In particular, for $\overline{X}_{1}$ and $\overline{X}_{s}$ the infinity is represented by $\{v = 0\}$ and $\{u = 0\}$, respectively.

\subsection{The expression of the compactified vector field}

The proof is carried out using a similar strategy as done in Section \ref{sec-one-useful-segment}, but performing compactification adapted to $\mathcal{P}$ instead of Poincaré--Lyapunov compactification. Indeed, we study the singularities in the set $\{uv = 0\}$ of the vector field $\overline{X}_{j}$ in an arbitrary chart $U_{j}$, for $j = 1,\dots, s$. Near each singularity at infinity, one must say under which conditions $\overline{X}_{j}$ and $\overline{(X_{G}^{U})}_{j}$ are topologically equivalent.

From the definition of compactification adapted to $\mathcal{P}$, there are three cases to consider when we deal with the change of coordinates induced by two consecutive vectors $\xi_{j-1}$, $\xi_{j}$ of the simple fan, for $j = 1,\dots,s$.
\begin{itemize}
    \item[(1)] The vectors $\xi_{j-1}$ and $\xi_{j}$ are not normal to no compact segment of $\mathcal{P}^{U}$.
    \item[(2)] The vector $\xi_{j-1}$ is not normal to any compact segment of $\mathcal{P}^{U}$, and $\xi_{j}$ is normal to a compact segment $\gamma_{k}^{U}$ of $\mathcal{P}^{U}$, for some $k = 1,\dots, l+1$.
    \item[(3)] The vector $\xi_{j-1}$ is normal to a compact segment of $\gamma_{k}^{U}$ of $\mathcal{P}^{U}$, for some $k = 0,\dots, l$, but $\xi_{j}$ is not normal to no compact segment of $\mathcal{P}^{U}$.
\end{itemize}

From the proof of \cite[Theorem A]{DOP24}, in case (1) the vector fields $\overline{X}_{j}$ and $\overline{(X_{G}^{U})}_{j}$ are topologically equivalent in the chart $U_{j}$, without requiring additional assumptions. Indeed, in this case, the only singularity at $\{uv = 0\}$ is the origin, and its topological behavior depends on the monomials related to a vertex of the upper diagram $\mathcal{P}^{U}$.

Cases (2) and (3) listed above can be treated in a similar fashion, so we carry out the computations for case (2). In this case, the vector $\xi_{j}$ is normal to a compact segment $\gamma_{k}\subset\mathcal{P}^{U}$ and, if the vector field $X_{\gamma_{k}^{U}}$ as defined in Equation \eqref{eq-quasi-upper} does not have singularities in $(\mathbb{R}^{*})^{2}$, then it follows from the proof of \cite[Theorem A]{DOP24} that in this chart $U_{j}$ the dynamics (under topological equivalence) is determined by the monomials related to points in the compact segment $\gamma_{k}^{U}\subset\mathcal{P}^{U}$. However, if $X_{\gamma_{k}^{U}}$ has singularities in $(\mathbb{R}^{*})^{2}$, then it might be necessary to consider monomials related to points below $\mathcal{P}^{U}$ as we shall see.

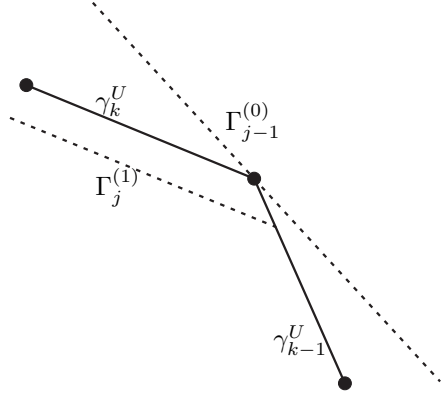
\begin{figure}[ht]\center{
\begin{overpic}[width=0.35\textwidth]{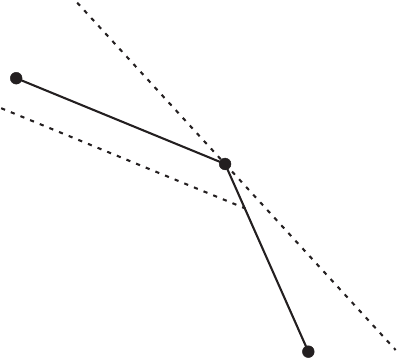}
%\begin{overpic}[grid,tics=10,width=0.35\textwidth]{fig-proof-general.pdf}

\put(20,65){\footnotesize{$\gamma_{k}^{U}$}}
\put(61,10){\footnotesize{$\gamma_{k-1}^{U}$}}
\put(50,60){\footnotesize{$\Gamma_{j-1}^{(0)}$}}
\put(20,45){\footnotesize{$\Gamma_{j}^{(1)}$}}

\end{overpic}}
\caption{\footnotesize{In the proof of Theorem \ref{mthm}, we consider auxiliary sets of the form $\Gamma_{j}^{(i)}$ defined in a similar fashion as in Section \ref{sec-def-gupp}. However, instead of constructing them from the collection of inward normal vectors $\Sigma^{*}_{\mathcal{P}}$, we construct them from the vectors of the simple fan $\Sigma_{\mathcal{P}}$. In this figure, we sketch the case where $\xi_{j}$ is the inward normal vector of some compact segment $\gamma_{k}^{U}$, whereas $\xi_{j-1}$ is normal to no edge.}}
\label{fig-proof-general}
\end{figure}

Recall that, in Section \ref{sec-def-gupp}, the numbers $\delta_{k}^{(i)}$ and the sets $\Gamma_{k}^{(i)}$ are defined considering the collection of inward normal vectors $\Sigma^{*}_{\mathcal{P}} = \{(\alpha_{i},\beta_{i})\}_{i=0}^{l+1}$. In what follows, however, we consider the numbers $\delta_{k}^{(i)}$ and the sets $\Gamma_{k}^{(i)}$ defined in a similar fashion as in Section \ref{sec-def-gupp}, however, the only difference is that, instead of constructing them from the collection of inward normal vectors $\Sigma^{*}_{\mathcal{P}} = \{(\alpha_{i},\beta_{i})\}_{i=0}^{l+1}$, we construct them from the vectors of the simple fan $\Sigma_{\mathcal{P}} = \{\xi_{j}\}_{j = 0}^{s}$. In this sense, since we are considering case (2) above, the set $\Gamma_{j}^{(0)}$ is contained in a segment $\gamma_{k}^{U}\subset\mathcal{P}^{U}$ for some $k$ (because $\xi_{j} = (\mu_{j},\nu_{j})$ is normal to such segment), whereas $\Gamma_{j-1}^{(0)}$ is only one vertex of $\mathcal{P}^{U}$ (because $\xi_{j-1}  = (\mu_{j-1},\nu_{j-1})$ is not normal to no segment). See Figure \ref{fig-proof-general}.

Fix an index $j$ such that case (2) holds. Firstly, we introduce some notation. It is important to remark that, in the computations below, the index $j$ is fixed, because it refers to the chart $U_{j}$ of the compactification.

Let $N_{j}$ be the smallest positive integer such that
\begin{equation*}
\mathcal{S}_{X} = \bigcup_{i = 0}^{N_{j}}\Gamma_{j}^{(i)},   
\end{equation*}
and write the polynomial vector field $X$ as
\begin{equation*}
X(x,y) = \displaystyle\sum_{i = 0}^{N_{j}}\left(P_{\Gamma_{j}^{(i)}}(x,y)\frac{\partial}{\partial x} + Q_{\Gamma_{j}^{(i)}}(x,y)\frac{\partial}{\partial y}  \right),    
\end{equation*}
where the monomials of $P_{\Gamma_{j}^{(i)}}$ and $Q_{\Gamma_{j}^{(i)}}$ are related to points of $\Gamma_{j}^{(i)}$. One further denote
\begin{equation*}
 P_{\Gamma_{j}^{(i)}}(x,y) = \displaystyle\sum_{\mu_{j}m + \nu_{j}n = \delta_{j}^{(i)}}a_{m,n}^{i,j}x^{m+1}y^{n}, \qquad Q_{\Gamma_{j}^{(i)}}(x,y) = \displaystyle\sum_{\mu_{j}m + \nu_{j}n = \delta_{j}^{(i)}}b_{m,n}^{i,j}x^{m}y^{n+1},
\end{equation*}
and recall that 
\begin{equation*}
   \delta_{j}^{(i)}  =  \displaystyle\min_{(m,n)\in\mathcal{S}\backslash \bigcup\Gamma_{j}^{(i-1)}}\{\alpha_{j} m + \beta_{j} n \}, \ \quad \
   \Gamma_{j}^{(i)}  =  \{(m,n)\in\mathcal{S}\backslash \cup\Gamma_{j}^{(i-1)} \ ; \ \alpha_{j} m + \beta_{j} n = \delta_{j}^{(i)}\}.   
\end{equation*}

After compactification and multiplication by $u^{|\delta_{j-1}^{(0)}|}v^{|\delta_{j}^{(0)}|}$, one obtains the vector field
\begin{equation}\label{eq-proof-compact-general}
\overline{X}_{j}(u,v) =   \displaystyle\sum_{i = 0}^{N_{j}}v^{|\delta_{j}^{(0)}| + \delta_{j}^{(i)}}\left( \displaystyle\sum_{(m,n)\in\Gamma_{j}^{(i)}}u^{|\delta_{j-1}^{(0)}| + \langle\xi_{j-1},(m,n)\rangle} \left( A_{m,n}^{i,j}u\frac{\partial}{\partial u} + B_{m,n}^{i,j}v\frac{\partial}{\partial v} \right)\right),  
\end{equation}
with
\begin{equation*}
A_{m,n}^{i,j} = \left(\nu_{j}a_{m,n}^{i,j} - \mu_{j}b_{m,n}^{i,j}\right), \qquad B_{m,n}^{i,j} = \left(-\nu_{j-1}a_{m,n}^{i,j} + \mu_{j-1}b_{m,n}^{i,j}\right).     
\end{equation*}

In Equation \eqref{eq-proof-compact-general}, recall that the index $j$ is fixed because it is related to the chart $U_{j}$. Powers of $v$ are indexed on $i$, which is related to the set $\Gamma_{j}^{(i)}$, and powers of $u$ depend on $(m,n)\in\Gamma_{j}^{(i)}$. As observed in the proof of \cite[Theorem A]{DOP24}, the singularities at infinity $\{uv = 0\}$ of the vector field $\overline{X}_{j}$ given in \eqref{eq-proof-compact-general} are, in fact, contained in the set $\{v = 0\}$. Indeed, this happens because the equality $\langle\xi_{j-1},(m,n)\rangle = \delta_{j - 1}^{(0)}$ holds only in the vertex $p_{k-1} = (m_{k-1},n_{k-1})$ of $\mathcal{P}^{U}$, which is exactly given by $\gamma_{k}^{U}\cap\gamma_{k-1}^{U}$, or also by $\Gamma_{j-1}^{(0)}\cap\Gamma_{j}^{(0)}$.

Let us study the behavior of the singularity at the origin of the chart $U_{j}$. Since $(a_{m_{k-1},n_{k-1}}^{0,j})^{2} + (b_{m_{k-1},n_{k-1}}^{0,j})^{2} \neq 0$ (because $p_{k-1}$ is a vertex) and $\operatorname{det}(\xi_{j-1},\xi_{j}) = 1$, then $(A_{m_{k-1},n_{k-1}}^{0,j})^{2} + (B_{m_{k-1},n_{k-1}}^{0,j})^{2} \neq 0$. Actually, one always has $B_{m_{k-1},n_{k-1}}^{0,j} \neq 0$, because if $B_{m_{k-1},n_{k-1}}^{0,j} = 0$ then $\{u = 0\}$ would be an arc of singularities, and we are considering the case where this does not happen (that is, all singularities at infinity are isolated). This means that the origin is either a hyperbolic singularity, or it is semi-hyperbolic whose center manifolds are tangent to $\{v = 0\}$, and the (un)stable manifold being transversal to $\{v = 0\}$. In both cases, the behavior depends on monomials related to $p_{k-1}\in\mathcal{P}^{U}$.

%both $A_{m_{k},n_{k}}^{0,j}$ and $B_{m_{k},n_{k}}^{0,j}$ are non zero, otherwise one would have an arc of singularities at infinity,  Therefore, it follows that the origin is a hyperbolic singularity, and then (due to the Grobman--Hartman Theorem) $\overline{X}_{j}$ is topologically equivalent to 
%\begin{equation*}
%A_{m_{k},n_{k}}^{0,j}u\frac{\partial}{\partial u} + B_{m_{k},n_{k}}^{0,j}v\frac{\partial}{\partial v}.   
%\end{equation*}

The next step is to study the singularities $(\lambda,0)\in\{v = 0\}$ with $\lambda \neq 0$, where $\lambda$ is a root of the polynomial equation
\begin{equation}\label{eq-root-general-case}
\displaystyle\sum_{(m,n)\in\Gamma_{j}^{(0)}}A_{m,n}^{0,j}u^{|\delta_{j-1}^{(0)}| + \langle\xi_{j-1},(m,n)\rangle + 1} = 0.     
\end{equation}

Due to item (3) of Definition \ref{def-mgupp-nd}, we must suppose that the polynomial in the left hand side of Equation \eqref{eq-root-general-case} is not identically zero, otherwise one would have arcs of singularities contained at infinity. In what follows, we further introduce the notation for the exponents of $u$ appearing in Equation \eqref{eq-proof-compact-general}:
\begin{equation*}
R_{m,n}^{i,j-1} := |\delta_{j-1}^{(i)}| + \langle\xi_{j-1},(m,n)\rangle.
\end{equation*}

If $\lambda\in\mathbb{R}$ is a non zero root of \eqref{eq-root-general-case}, we perform in Equation \eqref{eq-proof-compact-general} a change of coordinates of the form
\begin{equation*}
u = z + \lambda, \quad v = v,   
\end{equation*}
and then we get
\begin{equation}\label{eq-proof-compact-general-trans}
\overline{X}_{j,\lambda}(z,v)  =   \displaystyle\sum_{i = 0}^{N_{j}}\sum_{(m,n)\in\Gamma_{j}^{(i)}}v^{|\delta_{j}^{(0)}| + \delta_{j}^{(i)}} \left(  \displaystyle\sum_{r = 0}^{R_{m,n}^{i,j-1} + 1}z^{r}C_{m,n,r}^{i,j}\frac{\partial}{\partial z}  +  \displaystyle\sum_{r = 0}^{R_{m,n}^{i,j-1} }z^{r}D_{m,n,r}^{i,j}v\frac{\partial}{\partial v} \right),
\end{equation}
in which the coefficients $C_{m,n,r}^{i,j}$ and $D_{m,n,r}^{i,j}$ have the form
\begin{equation*}
C_{m,n,r}^{i,j} = \binom{R_{m,n}^{i,j-1} + 1}{r}A_{m,n}^{i,j}\lambda^{R_{m,n}^{i,j-1}  + 1 - r}, \qquad  D_{m,n,r}^{i,j} = \binom{R_{m,n}^{i,j-1}}{r}B_{m,n}^{i,j}\lambda^{R_{m,n}^{i,j-1} - r}.  
\end{equation*}

Observe that in Equation \eqref{eq-proof-compact-general-trans}, recall that the index $j$ is fixed because it is related to the chart $U_{j}$. Powers of $v$ are indexed on $i$, which is related to the set $\Gamma_{j}^{(i)}$. On the other hand, powers of $z$ are indexed on $r$, although $r$ depends on the indices $m,n,i$ (that is, $r:= r_{m,n}^{i,j-1}$).

We remark that several monomials associated with different points $(m,n)\in\Gamma_{j}^{(i)}$ may contribute to the same power of $z^{r}v^{|\delta_{j}^{(0)}| + \delta_{j}^{(i)}}$ the vector field \eqref{eq-proof-compact-general-trans}. In this sense, the coefficient of the monomial $z^{r}v^{|\delta_{j}^{(0)}| + \delta_{j}^{(i)}}$ (resp. $z^{r}v^{|\delta_{j}^{(0)}| + \delta_{j}^{(i)} + 1}$), as we vary $(m,n)\in\Gamma_{j}^{(i)}$, is the sum of coefficients of the form $C_{m,n,r}^{i,j}$ (resp. $D_{m,n,r}^{i,j}$). 

Finally, observe that setting $i = r = 0$ one obtains
\begin{equation*}
\displaystyle\sum_{(m,n)\in\Gamma_{j}^{(0)}}C_{m,n,0}^{0,j} = \displaystyle\sum_{(m,n)\in\Gamma_{j}^{(0)}}A_{m,n}^{0,j}\lambda^{R_{m,n}^{0,j-1}  + 1}  = 0,
\end{equation*}
because $\lambda$ satisfies Equation \eqref{eq-root-general-case}.

\subsection{Computing the topological normal form}

The next step is to apply Lemma \ref{prop-normal-form-semihyp}, following the same reasoning presented in Section \ref{sec-nf-one-seg}. Firstly, one must identify in the first component of \eqref{eq-proof-compact-general-trans} the coefficients of the linear term in the $z$ variable, and the term of lowest degree depending only on the $v$ variable. On the other hand, in the second component we pay attention to the monomial of the form $zv$ and the term of smallest degree depending only on the $v$ variable. 

Firstly, define the constants

\begin{equation*}
J_{m,n,1}^{0,j} = \sum_{(m,n)\in\Gamma_{j}^{(0)}}C_{m,n,1}^{0,j}, \qquad J_{m,n,0}^{1,j} = \sum_{(m,n)\in\Gamma_{j}^{(1)}}C_{m,n,0}^{1,j},
\end{equation*}
\begin{equation*}
L_{m,n,0}^{0,j} = \sum_{(m,n)\in\Gamma_{j}^{(0)}}D_{m,n,0}^{0,j}, \qquad L_{m,n,1}^{0,j} = \sum_{(m,n)\in\Gamma_{j}^{(0)}}D_{m,n,1}^{0,j}, \qquad L_{m,n,0}^{1,j} = \sum_{(m,n)\in\Gamma_{j}^{(1)}}D_{m,n,0}^{1,j}.
\end{equation*}

In this sense, in the first component we look at the monomials
\begin{equation*}
J_{m,n,1}^{0,j}z, \qquad \text{and} \qquad J_{m,n,0}^{1,j}v^{|\delta_{j}^{(0)}| + \delta_{j}^{(1)}},
\end{equation*}
and in the second component we look at the monomials 
\begin{equation*}
L_{m,n,0}^{0,j}v, \quad L_{m,n,1}^{0,j}zv, \quad \text{and} \quad L_{m,n,0}^{1,j}v^{|\delta_{j}^{(0)}| + \delta_{j}^{(1)} + 1}.
\end{equation*}

It is important to remark that the coefficients $J_{m,n,1}^{0,j}$, $L_{m,n,0}^{0,j}$ and  $L_{m,n,1}^{0,j}$ depend on monomials related to points of $\Gamma_{j}^{(0)}$, and the coefficients $J_{m,n,0}^{1,j}$ and $L_{m,n,0}^{1,j}$ depend on monomials related to points of $\Gamma_{j}^{(1)}$.

Now we proceed as in Section \ref{sec-nf-one-seg}. If $L_{m,n,0}^{0,j} \neq 0$, then the singularity $\lambda$ is either hyperbolic or semi-hyperbolic. If it is hyperbolic, we are done by the Grobman--Hartman Theorem. If it is semi-hyperbolic, central manifolds are tangent to infinity, and the (un)stable manifold is transversal to infinity. In both cases, the behavior of the vector field near the singularity $\lambda$ depends only on the monomials associated to the set $\Gamma_{j}^{(0)}$.

\begin{remark}
By the proof of \cite[Theorem A]{DOP24}, if $X_{\gamma_{k}^{U}}$ does not have singularities in $(\mathbb{R}^{*})^{2}$, then this implies that $L_{m,n,0}^{0,j} \neq 0$. 
\end{remark}

However, we allow the case where $X_{\gamma_{k}^{U}}$ has singularities in $(\mathbb{R}^{*})^{2}$. So, if $L_{m,n,0}^{0,j} = 0$ then one must require $J_{m,n,1}^{0,j} \neq 0$, that is, the singularity is semi-hyperbolic and the (un)stable manifold is tangent to infinity. In this case, it follows by Lemma \ref{prop-normal-form-semihyp} that if there exist a polynomial function $F_{\lambda,j}$ such that the Vector Field \eqref{eq-proof-compact-general-trans} is topologically equivalent to
\begin{equation*}
J_{m,n,1}^{0,j}z\frac{\partial}{\partial z} + F_{\lambda,j}\left(L_{m,n,0}^{1,j}, J_{m,n,0}^{1,j},L_{m,n,1}^{0,j}\right) v^{|\delta_{j}^{(0)}| + \delta_{j}^{(1)} + 1}\frac{\partial}{\partial v},  
\end{equation*}
provided that $F_{\lambda,j}\left(L_{m,n,0}^{1,j}, J_{m,n,0}^{1,j},L_{m,n,1}^{0,j}\right) \neq 0$. Recall that the coefficients $J_{m,n,1}^{0,j}$ and $L_{m,n,1}^{0,j}$ depend on $\Gamma_{j}^{(0)}$, and the coefficients $J_{m,n,0}^{1,j}$ and $L_{m,n,0}^{1,j}$ depend on $\Gamma_{j}^{(1)}$.

Finally, for each singularity $(\lambda,0)$ at infinity of the chart $U_{j}$, define the set
\begin{equation}\label{eq-def-u1}
 \mathfrak{U}_{2,\lambda}^{j}(\mathcal{P}) := \left\{ X\in\mathfrak{U}_{1}(\mathcal{P}); \ L_{m,n,0}^{0,j}= 0 \right\} \cap \left\{X\in\mathfrak{U}_{1}(\mathcal{P}); \ F_{\lambda,j}\left(L_{m,n,0}^{1,j}, J_{m,n,0}^{1,j},L_{m,n,1}^{0,j}\right) = 0 \right\}.   
\end{equation}

It is clear from our reasoning that if $X\not\in \mathfrak{U}_{2,\lambda}^{j}(\mathcal{P})$, then the compactifications of $X$ and $X_{G}^{U}$ in the chart $U_{j}$ are topologically equivalent near $(\lambda,0)$. We further define
\begin{equation}\label{eq-def-u1-final}
 \mathfrak{U}_{2}^{j}(\mathcal{P}) := \bigcup_{\lambda\in\mathcal{Z}_{j}} \mathfrak{U}_{2,\lambda}^{j}(\mathcal{P}),  \quad \mathfrak{U}_{2}(\mathcal{P}) := \bigcup_{j = 0}^{l+1}\mathfrak{U}_{2}^{j}(\mathcal{P}), \quad \widetilde{\mathfrak{U}}_{1}(\mathcal{P}) := \mathfrak{U}_{1}(\mathcal{P})\backslash \mathfrak{U}_{2}(\mathcal{P}),
\end{equation}
where $\mathcal{Z}_{j}$ is the set of all real roots of Equation \eqref{eq-root-general-case} in this chart $U_{j}$ (that is, $\lambda$ is a singularity at infinity of the chart $U_{j}$, for some $j = 0,\dots, l+1$). Thus, if $X\in\widetilde{\mathfrak{U}}_{1}(\mathcal{P})$, then the compactifications of $X$ and $X_{G}^{U}$ are topologically equivalent near each singularity $(\lambda,0)$ at infinity. The set $\widetilde{\mathfrak{U}}_{1}(\mathcal{P})$ is open and dense in $\mathfrak{U}_{1}(\mathcal{P})$ with respect to the Zariski topology. Finally, observe that the subsets $\widetilde{\mathfrak{U}}_{1}(\mathcal{P})$ and $\mathfrak{U}_{2}(\mathcal{P})$ of $\mathfrak{U}_{1}(\mathcal{P})$ do not depend on the initial vector field $X$ given, because each set $\mathfrak{U}_{2,\lambda}^{j}(\mathcal{P})$ is defined independently of $X$.

\begin{remark}
Observe that, in the case $L_{m,n,0}^{0,j} = J_{m,n,1}^{0,j} = 0$, the singularity $(\lambda,0)\in\{v = 0\}$ is neither hyperbolic nor semi-hyperbolic, therefore it is non-elementary. This would lead to further blow-up analysis and we do not treat this case in this paper.  
\end{remark}

The existence of the homeomorphism that gives the topological equivalence in a neighborhood of infinity follows with a completely analogous reasoning as described in Section \ref{sec-existence-top-eq}. This completes the proof.

\section{Concluding remarks}\label{sec-concluding-remarks}

The goal of this paper is to generalize the results of  \cite{DOP24, OliVal26} concerning topological normal forms for the study of the dynamics at infinity, and our approach is inspired in \cite{Zup96,Zup00}. More precisely, the topological normal form addressed in this paper is the \textit{minimal generalized upper principal part}. In order to define such a topological normal form, firstly, in Section \ref{sec-def-gupp}, we described the so called \textit{Newton decomposition} (which is, in general, different from the classical \textit{Newton filtration}). This allows us to write $X$ as the sum of polynomial vector fields
\begin{equation*}
X(x,y) = \displaystyle\sum_{i = 0}^{N}X_{i}^{U}(x,y),    
\end{equation*}
where the monomials of $X_{i}^{U}$ are related to points contained in $\Gamma^{(i)}$ and $X_{0}^{U} = X_{\Delta}^{U}$, with $i = 0,\dots, N$.
Next, we defined the \textit{generalized upper principal part} $X_{\Gamma}^{U}$ of $X$ as
\begin{equation*}
X_{\Gamma}^{U} = X_{0}^{U} + X_{1}^{U}.    
\end{equation*}

Following Section \ref{sec-def-gupp}, we saw that some monomials of $X_{\Gamma}^{U}$ can be dropped as described in Definition \ref{def-mgupp}, which leads to the definition of \textit{minimal} generalized upper principal part $X^{U}_{G}$.

The main result of the paper is Theorem \ref{mthm}, and it is stated in Section \ref{sec-main-thms}. It says that, given a Newton polygon $\mathcal{P}$, there exist an open and dense subset $\widetilde{\mathfrak{U}}_{1}(\mathcal{P})\subset \mathfrak{U}_{1}(\mathcal{P})$ satisfying the following property: If $X\in \widetilde{\mathfrak{U}}_{1}(\mathcal{P})$ and its minimal generalized upper principal part $X_{G}^{U}$ is non-degenerate (in the sense of Definition \ref{def-mgupp-nd}), then $X$ and $X_{G}^{U}$ are topologically equivalent near infinity. 

The three conditions required in Definition \ref{def-mgupp-nd} can be checked performing \textit{Bendixson compactification} in the vector field $X$ given. More precisely, if $\mathcal{B}(X)$ is the compactified vector field according to Bendixson, then the origin is an isolated, non monodromic and non dicritical singularity of $\mathcal{B}(X)$. This was already observed in Section \ref{sec-main-thms}. Now,  specifically concerning the condition (3) of Definition \ref{def-mgupp-nd}, this condition can be algebraically checked. Indeed, it is equivalent to say that the polynomial in the left hand side of Equation \eqref{eq-root-general-case} is not identically zero, for each chart of the compactification procedure adopted in the proof presented in Section \ref{sec-proofs}. In the case where $\mathcal{P}^{U}$ has only one useful segment, such condition can be simplified even more: it is equivalent to say that the polynomial in the left hand side of Equation \eqref{eq-proof-roots} is not identically zero. Finally, we remark that the set $\widetilde{\mathfrak{U}}_{1}(\mathcal{P})\subset \mathfrak{U}_{1}(\mathcal{P})$ does not depend on the vector field $X$ given.

The proof of Theorem \ref{mthm} is given in Section \ref{sec-one-useful-segment} for the case where $\mathcal{P}$ has only one useful segment, and in Section \ref{sec-proofs} for the general case. In both cases, the idea of the proof relies in the same principle: study the singularities at infinity by applying Lemma \ref{prop-normal-form-semihyp}. Indeed, when a singularity at infinity is semi-hyperbolic whose central manifolds are transversal to infinity, one should check if the coefficients of the translated vector field satisfies the assumptions of Lemma \ref{prop-normal-form-semihyp}. This can be checked by simply using the polynomial $F$ defined in the statement of Lemma \ref{prop-normal-form-semihyp}. The polynomial $F$ is important when we define the subset $\widetilde{\mathfrak{U}}_{1}(\mathcal{P})\subset \mathfrak{U}_{1}(\mathcal{P})$ in Equations \eqref{eq-def-u1} and \eqref{eq-def-u1-final}. Finally, recall that in the proof given in Section \ref{sec-one-useful-segment} we performed Poincaré--Lyapunov compactification, whereas in the proof of the general case given in Section \ref{sec-proofs} we performed \textit{compactification adapted to $\mathcal{P}$} (which is toric compactification) instead.

Our first remark is that our approach suggests an alternative proof for the main result of \cite{Zup96,Zup00}. Indeed, instead of proceeding by induction and applying the algorithm described by Pelletier \cite{Pel95}, one can perform a \textit{toric resolution of singularities} as Brunella and Miari \cite{BM90}, and then applying the Normal Form Theorem to study the singularities along the exceptional divisor. We believe that the steps and the computations would be completely analogous as presented in this paper. In our paper, the compactification adapted to $\mathcal{P}$ was useful in order to precisely define the subset $\widetilde{\mathfrak{U}}_{1}(\mathcal{P})\subset \mathfrak{U}_{1}(\mathcal{P})$ in Equations \eqref{eq-def-u1} and \eqref{eq-def-u1-final}.

The second remark is that one might wonder if, for planar polynomial vector fields belonging to $\mathfrak{U}_{2}$, one could prove the topological equivalence at infinity between $X$ and the polynomial vector field
\begin{equation*}
X^{(2)}(x,y) = X_{0}^{U}(x,y) + X_{1}^{U}(x,y) + X_{2}^{U}(x,y),    
\end{equation*}
where $X_{2}^{U}$ has monomials related to points of $\Gamma^{(2)}$. Actually, one could ask a more general question. For $k = 1,\dots, N$, we inductively define $\mathfrak{U}_{k} = \widetilde{\mathfrak{U}}_{k}\cup\mathfrak{U}_{k+1}$, where $\widetilde{\mathfrak{U}}_{k}$ is open and dense in $\mathfrak{U}_{k}$ and, for $X\in\widetilde{\mathfrak{U}}_{k}$, the vector field $X$ is topologically equivalent near infinity to
\begin{equation*}
\begin{array}{rlc}
     X^{(k)}(x,y) & = & X_{0}^{U}(x,y) + X_{1}^{U}(x,y) + \dots + X_{k}^{U}(x,y) \\
     & = & X_{\Delta}^{U}(x,y) + X_{1}^{U}(x,y) + \dots + X_{k}^{U}(x,y), 
\end{array}  
\end{equation*}
where $X_{i}^{U}$ has monomials related to points of $\Gamma^{(i)}$, for $i = 0,\dots, k$. If we use the approach presented here, the proof would be more involved and delicate. Indeed, Lemma \ref{prop-normal-form-semihyp} is very important in the proof of our Theorem \ref{mthm}, and in order to prove such conjecture it would be necessary to prove an analogous result to Lemma \ref{prop-normal-form-semihyp} but considering higher order terms of the analytic vector field given in Equation \eqref{eq-vf-prop}. This is a topic for future study.

%%%%%%%%%%%%%%%%%%%%%%%%%%%%%%%%%%%%%%%%%%%%
%%%%%%%%%%%%%%%%%%%%%%%%%%%%%%%%%%%%%%%%%%%%
%%%%%%%%%%%%%%%%%%%%%%%%%%%%%%%%%%%%%%%%%%%%
%%%%%%%%%%%%%%%%%%%%%%%%%%%%%%%%%%%%%%%%%%%%

\section*{Declarations}
\textbf{Conflict of interest} The authors declare that they have no conflict of interest.

\textbf{Data Availability Statement} Data sharing not applicable to this article as no datasets were generated or
analysed during the current study.

\section*{Acknowledgments and funding}

The authors thank Vesna \v{Z}upanovi\'c for reading a preliminary version of this manuscript and for her valuable comments.

T.M. Dalbelo is supported by FAPESP grants 2019/21181-0 and 2024/22060-0 and by CNPq grant 403959/2023-3. R. Oliveira was financed, in part, by the São Paulo Research Foundation (FAPESP), Brazil, process numbers 2019/21181-0; the Brazilian CNPq grant number 310857/2023-6 and 407454/2023-3. O.H. Perez is supported by Sao Paulo Research Foundation (FAPESP) grant 2021/10198-9.

\appendix

\section{The Normal Form Theorem}\label{appendix-nft}

In what follows it is presented a brief introduction on the Normal Form Theorem. The proof can be found in Takens' paper \cite{Tak74} and in Bruno's book \cite{Bru89} (see Sections 1.5 and 1.9 of Chapter II in this last reference). A detailed exposition of the Normal Form Theorem can also be found in the textbooks \cite{AP90, DLA06, GH83}. For the sake of generality, we present the Normal Form Theorem (Theorem \ref{teo-nft}) in the case that $X$ is a vector field defined in $\mathbb{R}^{n}$.

\subsection{Statement of the Normal Form Theorem}\label{sec-nft}

Consider a $C^{r}$ vector field (with $r\geq 1$) defined in a neighborhood $U$ of $0\in\mathbb{R}^{n}$ of the form
\begin{equation*}
X(\textbf{x}) = A\textbf{x} + X_{2}(\textbf{x}) + \dots + X_{r}(\textbf{x}) + \left(||\textbf{x}||^{r+1}\right),
\end{equation*}
in which $X_{i}(\textbf{x})$ is a homogeneous vector field of degree $i = 2,\dots, r$ and $A$ is a non-zero real valued matrix. We define the \textit{adjoint operator} (or \textit{Lie bracket operator}, or \textit{homological operator}) as
$$
\begin{array}{rcl}
   L_{i}A: \mathcal{H}_{i}\left(\mathbb{R}^{n}\right)  & \rightarrow & \mathcal{H}_{i}\left(\mathbb{R}^{n}\right) \\
    Y & \mapsto & [A\textbf{x},Y] = J_{Y}A\textbf{x} - J_{A\textbf{x}}Y
\end{array}
$$
in which $\mathcal{H}_{i}\left(\mathbb{R}^{n}\right)$ is the space of all homogeneous vector fields of degree $i$, and $J_{Y}$ is the Jacobian matrix of $Y$. We further denote $\mathcal{R}_{i} = L_{i}A\left(\mathcal{H}_{i}\left(\mathbb{R}^{n}\right)\right)$ and $\mathcal{G}_{i}$ is a set such that $\mathcal{H}_{i}\left(\mathbb{R}^{n}\right) = \mathcal{R}_{i}\oplus \mathcal{G}_{i}$.

\begin{theorem}[Normal Form Theorem]\label{teo-nft} Let $X$ be a $C^{r}$ vector field defined in a neighborhood $U$ of $0\in\mathbb{R}^{n}$, with $X(0)=0$ and $J_{X}(0) = A\neq 0$. Let $\mathcal{R}_{i}$ and $\mathcal{G}_{i}$ be as above. Then there is an analytic change of coordinates $\Psi:U\rightarrow \mathbb{R}^{n}$ such that $Y = \Psi^{*}X$ is of the form
\begin{equation*}
Y(\textbf{y}) = A\textbf{y} + Y_{2}(\textbf{y}) + \dots + Y_{r}(\textbf{y}) + O\left(||\textbf{y}||^{r+1}\right),
\end{equation*}
 with $Y_{i}\in \mathcal{G}_{i}$ for $i = 2,\dots, r$.
\end{theorem}

The proof of Theorem \ref{teo-nft} gives us an algorithm for computing normal forms of $C^{r}$ vector fields. This is an iterative process. At each step, we perform a change of coordinates of the form
$$\Psi_{i}(\textbf{x}) = \textbf{x} + \psi_{i}(\textbf{x}),$$
that is, the diffeomorphism has the form ``identity plus homogeneous polynomial of degree $i$''. The vector field $\psi_{i}(\textbf{x})\in \mathcal{H}_{i}\left(\mathbb{R}^{n}\right)$ must be well-chosen in the sense that we want to eliminate as many monomials as possible. Terms that cannot be eliminated by this method are called \textit{resonant terms}. A vector field whose Taylor series contains only resonant terms is called \textit{normal form}.

Applying this algorithm in the case of a semi-hyperbolic singularity of planar vector field, and assuming without loss of generality that $A$ is in its Jordan form, one obtains the following result (see also \cite{AP90, DLA06, GH83}).

\begin{proposition}\label{prop-analytic-nf}
If the planar vector field $X$ given in Equation \eqref{eq-vf-prop} has an isolated semi-hyperbolic singularity positioned at the origin, and assuming without loss of generality that the linearization of $X$ is in its Jordan form, then there exists an analytic change of coordinates $\Phi$ such that $Y = \Phi^{*}X$, where form some integer $m\geq 2$ it holds
\begin{equation}\label{eq-analytic-nf}
Y(x,y) = \left(a_{1,0}x + \displaystyle\sum_{s = 1}^{m-1}\widetilde{a}_{1,s}xy^{s}\right)\frac{\partial}{\partial x} + \widetilde{b}_{0,m}y^{m}\frac{\partial}{\partial y} + \left( ||(x,y)||^{m+1}\right)\left(\frac{\partial}{\partial x} + \frac{\partial}{\partial y}\right).   
\end{equation}
\end{proposition}

\subsection{Proof of Lemma \ref{prop-normal-form-semihyp}}\label{proof-prop}

We are interested in topological normal forms of semi-hyperbolic singularities of planar vector fields of the form
\begin{equation*}
X(x,y) = X_{1}(x,y)\frac{\partial}{\partial x} + X_{2}(x,y)\frac{\partial}{\partial y} = \left(\sum_{r + s \geq 1}a_{r,s}x^{r}y^{s}\right)\frac{\partial}{\partial x} + \left(\sum_{r + s \geq 1}b_{r,s}x^{r}y^{s}\right)\frac{\partial}{\partial y},
\end{equation*}
as given in Equation \eqref{eq-vf-prop}. It follows from  \cite[Chapter 2]{DLA06} that the topological behavior of the semi-hyperbolic singularity is given by a vector field of the form
\begin{equation*}
\widetilde{Y}(x,y) = a_{1,0}x\frac{\partial}{\partial x} + by^{m}\frac{\partial}{\partial y}.
\end{equation*}

Intuitively, the topological behavior is determined by the linear term in the first component, and on the term of smallest degree depending only on $y$ in the second component. Moreover, depending on the sign of $a_{1,0},b$ and the parity of $m$, the singularity can be a saddle, a node or a saddle-node. In order to find such a topological normal form, firstly we apply the algorithm described in Section \ref{sec-nft}.

Suppose that $l\geq 1$ is the smallest integer such that $X$ given in \eqref{eq-first-vf} has the monomial $a_{0,l}y^{l}\frac{\partial}{\partial x}$ or $b_{0,l+1}y^{l+1}\frac{\partial}{\partial y}$. We divide the proof in two cases.

\

\noindent\textbf{Case $l = 1$ and $b_{2,0} = 0$.} The vector field $X$ is given by 
\begin{equation}\label{eq-first-vf-pt2}
X(x,y) = \left(a_{1,0}x + a_{0,1}y + \sum_{r + s \geq 2}a_{r,s}x^{r}y^{s}\right)\frac{\partial}{\partial x} + \left(b_{0,2}y^{2} + b_{1,1}xy +  \sum_{r + s \geq 3}b_{r,s}x^{r}y^{s}\right)\frac{\partial}{\partial y}.
\end{equation}

In this case, the (un)stable manifold is tangent to $\{y = 0\}$, whereas center manifolds are not. We perform a linear change of coordinates $(a_{1,0}x + a_{0,1}y, y) \mapsto (\tilde{x},\tilde{y})$, and we obtain the vector field
\begin{equation*}
\begin{array}{rcl}
     \widetilde{X}(\tilde{x},\tilde{y}) & = & \left(
  \begin{array}{cc}
    a_{1,0} & a_{0,1} \\
    0 & 1 \\
  \end{array}
\right)X\left(\frac{\tilde{x} - a_{0,1}\tilde{y}}{a_{1,0}},\tilde{y}\right) \\
     & = & \left(a_{1,0}\tilde{x} + \displaystyle\sum_{r + s \geq 2}\tilde{a}_{r,s}\tilde{x}^{r}\tilde{y}^{s}\right)\displaystyle\frac{\partial}{\partial \tilde{x}} + \left( \left(b_{0,2} - \frac{b_{1,1}a_{0,1}}{a_{1,0}}\right)\tilde{y}^{2} + b_{1,1}\tilde{x}\tilde{y} +  \displaystyle\sum_{r + s \geq 3}\tilde{b}_{r,s}\tilde{x}^{r}\tilde{y}^{s}\right)\displaystyle\frac{\partial}{\partial \tilde{y}}, 
\end{array}
\end{equation*}
for some new coefficients $\tilde{a}_{r,s},\tilde{b}_{r,s}$. Assuming $b_{0,2} - \frac{b_{1,1}a_{0,1}}{a_{1,0}} \neq 0$, we apply Proposition \ref{prop-analytic-nf}, and therefore we obtain the topological normal form given in the first item of Lemma \ref{prop-normal-form-semihyp}.

\

\noindent\textbf{Case $l\geq 2$.} Here, it is not necessary to suppose $b_{2,0} = 0$. Since the origin is a semi-hyperbolic singularity, we write it as
\begin{equation}\label{eq-first-vf}
\begin{array}{rcrcl}
    X(x,y) & = & X_{1}(x,y)\displaystyle\frac{\partial}{\partial x} & + & X_{2}(x,y)\displaystyle\frac{\partial}{\partial y} \\
     & = & \left(a_{1,0}x + \displaystyle\sum_{r + s \geq 2}a_{r,s}x^{r}y^{s}\right)\displaystyle\frac{\partial}{\partial x} & + & \left(\displaystyle\sum_{r + s \geq 2}b_{r,s}x^{r}y^{s}\right)\displaystyle\frac{\partial}{\partial y}.
\end{array}
\end{equation}

We start the Normal Form algorithm described in Section \ref{sec-nft} with a degree of homogeneity $r'+s' = h \geq 2$. Then, one should find a suitable change of coordinates of the form
\begin{equation*}
\Psi_{h}(\tilde{x},\tilde{y}) = \left(\Psi_{1,h}(\tilde{x},\tilde{y}),\Psi_{2,h}(\tilde{x},\tilde{y})\right) = \left(\tilde{x} + \displaystyle\sum_{r'+s' = h}c_{r',s'}\tilde{x}^{r'}\tilde{y}^{s'}, \quad \tilde{y} + \displaystyle\sum_{r'+s' = h}d_{r',s'}\tilde{x}^{r'}\tilde{y}^{s'}\right).    
\end{equation*}

The notation $r'$ and $s'$ is adopted in order to highlight that such powers come from the diffeomorphism $\Psi_{h}$. By the Normal Form algorithm applied in the case of a semi-hyperbolic singularity, for each pair $(r',s')$ satisfying $r'+s' = h \geq 2$ the coefficients are given by
\begin{equation}\label{eq-proof-coef}
c_{r',s'} = \frac{a_{r',s'}}{a_{1,0}(r'-1)}, \quad d_{r',s'} = \frac{b_{r',s'}}{a_{1,0}r'}.     
\end{equation}

By the expression of the coefficients \eqref{eq-proof-coef}, it is clear that terms of the form $xy^{s'}$ in the first component cannot be eliminated by this algorithm, whereas terms of the form $y^{s'}$ in the second component also cannot be eliminated.

Denoting the Jacobian matrix of $\Psi_{h}$ by $J\Psi_{h}$, one obtains the vector field
\begin{equation*}
\begin{array}{rcl}
\widetilde{X}(\tilde{x},\tilde{y}) & = & \left(J\Psi_{h}(\tilde{x},\tilde{y})\right)^{-1}X\left(\Psi_{h}(\tilde{x},\tilde{y})\right) \\
& = & \displaystyle\frac{1}{\operatorname{det}J\Psi_{h}(\tilde{x},\tilde{y})}\left(
  \begin{array}{cc}
    \frac{\partial \Psi_{2,h}}{\partial y} & -\frac{\partial \Psi_{1,h}}{\partial y} \\
    -\frac{\partial \Psi_{2,h}}{\partial x} & \frac{\partial \Psi_{1,h}}{\partial x} \\
  \end{array}
\right)X\left(\Psi_{h}(\tilde{x},\tilde{y})\right).
\end{array}
\end{equation*}

Since $\operatorname{det}J\Psi_{h}(\tilde{x},\tilde{y}) = 1 + O(u,v)$, by performing a time reparameterization it follows that near the origin the vector field $\widetilde{X}$ is topologically equivalent to
\begin{equation}\label{eq-final-vf}
\widetilde{Y}(\tilde{x},\tilde{y}) =  \left(
  \begin{array}{cc}
    \frac{\partial \Psi_{2,h}}{\partial y} & -\frac{\partial \Psi_{1,h}}{\partial y} \\
    -\frac{\partial \Psi_{2,h}}{\partial x} & \frac{\partial \Psi_{1,h}}{\partial x} \\
  \end{array}
\right)X\left(\Psi_{h}(\tilde{x},\tilde{y})\right).
\end{equation}

The next step is to study the coefficient of the term of smallest degree depending only on $\tilde{y}$ in the second component of $\widetilde{Y}$, which is given by
\begin{equation}\label{eq-proof-second-component}
\left(-\frac{\partial \Psi_{2,h}}{\partial x}X_{1}\left(\Psi_{h}(\tilde{x},\tilde{y})\right) + \frac{\partial \Psi_{1,h}}{\partial x}X_{2}\left(\Psi_{h}(\tilde{x},\tilde{y})\right)\right)\frac{\partial}{\partial \tilde{y}}.  
\end{equation}

Firstly, we analyze the terms depending only on $y$ in the part $\frac{\partial \Psi_{1,h}}{\partial x}X_{2}\left(\Psi_{h}(\tilde{x},\tilde{y})\right)$ in Equation \eqref{eq-proof-second-component}. Since $\frac{\partial \Psi_{1,h}}{\partial x} = 1 + O(u,v)$, then 
\begin{equation*}
\frac{\partial \Psi_{1,h}}{\partial x}X_{2}\left(\Psi_{h}(\tilde{x},\tilde{y})\right) = \left(1 + O(u,v)\right)\left(\sum_{r + s \geq 2} b_{r,s}\left(\tilde{x} + \displaystyle\sum_{r'+s' = h}c_{r',s'}\tilde{x}^{r'}\tilde{y}^{s'}\right)^{r}\left(\tilde{y} + \displaystyle\sum_{r'+s' = h}d_{r',s'}\tilde{x}^{r'}\tilde{y}^{s'}\right)^{s} \right).   
\end{equation*}

Recall that the notation $r'$ and $s'$ is adopted in order to stress that such powers come from the diffeomorphism $\Psi_{h}$. Moreover, $r' + s' = h \geq 2$. Therefore, in this part the terms of smaller degree depending only on $\tilde{y}$ have the form $b_{r,s}(c_{0,h})^{r}\tilde{y}^{hr+s}$, with $r+s \geq 2$.

Let us check the terms depending only on $\tilde{y}$ in the part $-\frac{\partial \Psi_{2,h}}{\partial x}X_{1}\left(\Psi_{h}(\tilde{x},\tilde{y})\right)$ in Equation \eqref{eq-proof-second-component}. It follows that 
\begin{equation*}
\begin{array}{rcl}
     -\frac{\partial \Psi_{2,h}}{\partial x}X_{1}\left(\Psi_{h}(\tilde{x},\tilde{y})\right) & = & -\left(\displaystyle\sum_{r'+s' = h}d_{r',s'}r'\tilde{x}^{r'-1}\tilde{y}^{s'}\right)\left[a_{1,0}\left(\tilde{x} + \displaystyle\sum_{r'+s' = h}c_{r',s'}\tilde{x}^{r'}\tilde{y}^{s'}\right)\right. + \\
   & + & \left. \displaystyle\sum_{r + s \geq 2} a_{r,s}\left(\tilde{x} + \displaystyle\sum_{r'+s' = h}c_{r',s'}\tilde{x}^{r'}\tilde{y}^{s'}\right)^{r}\left(\tilde{y} + \displaystyle\sum_{r'+s' = h}d_{r',s'}\tilde{x}^{r'}\tilde{y}^{s'}\right)^{s}\right]
\end{array}
\end{equation*}

Then the terms of smaller degree depending only on $\tilde{y}$ in the first part of Equation \eqref{eq-proof-second-component} will have the form
\begin{equation*}
-d_{1,h-1}a_{1,0}c_{0,h}\tilde{y}^{2h-1}, \quad -d_{1,h-1} a_{r,s}(c_{0,h})^{r}\tilde{y}^{(h - 1) + hr+s}, \qquad r + s \geq 2.   
\end{equation*}

Since $h\geq 2$, then it is always true that $(h - 1) + hr+s > hr+s$. Gathering the information obtained in the first and second parts of Equation \eqref{eq-proof-second-component}, finally we obtain the candidates for term of smaller degree depending only on $\tilde{y}$, which are given by
\begin{equation*}
-d_{1,h-1}a_{1,0}c_{0,h}\tilde{y}^{2h-1}, \quad b_{r,s}(c_{0,h})^{r}\tilde{y}^{hr+s}, \qquad r + s \geq 2.       
\end{equation*}

Suppose that $l\geq 2$ is the smallest integer such that $X$ given in \eqref{eq-first-vf} has the monomial $a_{0,l}y^{l}\frac{\partial}{\partial x}$ or $b_{0,l+1}y^{l+1}\frac{\partial}{\partial y}$, and recall the expressions of the coefficients $c_{r',s'}, d_{r',s'}$ in Equation \eqref{eq-proof-coef}. Then, after performing the change of coordinates $\Psi_{h}$, the terms in the second component of \eqref{eq-final-vf} of smallest with respect to the $\tilde{y}$ variable have the form
\begin{equation*}
\begin{array}{cl}
    b_{1,h-1}\frac{a_{0,h}}{a_{1,0}}\tilde{y}^{2h-1}, &  \\
   \qquad b_{0,l+1}\tilde{y}^{l+1},  & \text{if $r = 0$ and $s = l+1$}, \\
   \qquad -b_{1,1}\frac{a_{0,l}}{a_{1,0}}\tilde{y}^{l+1}, & \text{if $r = 1$, $h = l$ and $s = 1$}.
\end{array}
\end{equation*}

Since $l$ is the smallest integer such that $a_{0,l} \neq 0$ or $b_{0,l+1} \neq 0$, the term $b_{1,h-1}\frac{a_{0,h}}{a_{1,0}}\tilde{y}^{2h-1}$ should be taken into account when $h = l$, which gives us the term $b_{1,l-1}\frac{a_{0,l}}{a_{1,0}}\tilde{y}^{2l-1}$. However, observe that $2l -1 \geq l+1$ and the equality holds for $l = 2$. 

If $l \geq 3$, the term of smallest degree depending only on $\widetilde{y}$ in the second component of $\widetilde{Y}$ given in Equation \eqref{eq-final-vf} is given by $\left(b_{0,l+1} - \frac{a_{0,l}b_{1,1}}{a_{1,0}}  \right)\tilde{y}^{l+1} \frac{\partial}{\partial \tilde{y}}$, provided that such coefficient is nonzero. When $l = 2$, then the term of smallest degree depending only on $\widetilde{y}$ in the second component of $\widetilde{Y}$ given in Equation \eqref{eq-final-vf} is given by $b_{0,3}\tilde{y}^{3} \frac{\partial}{\partial \tilde{y}}$.

In both cases, it follows by Proposition \ref{prop-analytic-nf} that such term of smallest degree in the second component cannot be eliminated by the algorithm described in Section \ref{sec-nft}. Finally, it follows from results of \cite[Chapter 2]{DLA06} that the topological normal form is given as in Lemma \ref{prop-normal-form-semihyp}. This completes the proof.

Of course, when the coefficients of the vector field of Lemma \ref{prop-normal-form-semihyp} satisfy $F\left(b_{0,l+1}, a_{0,l}, b_{1,1}\right) = 0$, then one must continue the Normal Form algorithm taking into account higher order terms of $X$. This would lead to taking into account my different sub-cases and then the computations will be more delicate. We leave this for future study.

\section{Proof of Lemma \ref{lemma-compact-x-trans}}\label{appendix-proofs}

It is convenient to introduce the following notation.
\begin{equation}\label{eq-expression-pq}
P_{0}^{U}(1,u) = \displaystyle\sum_{n = 0}^{\alpha + 1}a_{n}u^{n}, \quad P_{1}^{U}(1,u) = \displaystyle\sum_{n = 0}^{N_{1}}c_{n}u^{n}, \quad Q_{0}^{U}(1,u) = \displaystyle\sum_{n = 0}^{\alpha + 1}b_{n}u^{n}, \quad Q_{1}^{U}(1,u) = \displaystyle\sum_{n = 0}^{N_{2}}d_{n}u^{n},
\end{equation}

We remark that in Case I, we set $a_{\alpha} = a_{\alpha + 1} = b_{\alpha} = b_{\alpha + 1} = 0$, but it holds that $a_{\alpha - 1} \neq 0 \neq b_{0}$. Indeed, the component $P_{0}^{U}$ has a monomial of the form $a_{\alpha - 1}y^{\alpha -1}\frac{\partial}{\partial x}$, and the component $Q_{0}^{U}$ also has monomial of the form $b_{0}x^{\beta -1}\frac{\partial}{\partial y}$. In Case II, we set $a_{\alpha + 1} = b_{0} = 0$, and it holds that $a_{\alpha}^{2} + b_{\alpha + 1}^{2} \neq 0$ and $a_{0}^{2} + b_{1}^{2} \neq 0$. This is true because the component $P_{0}^{U}$ has monomials of the form $a_{\alpha}xy^{\alpha}\frac{\partial}{\partial x}$ and/or $a_{0}x^{\beta +1}\frac{\partial}{\partial x}$, and the component $Q_{0}^{U}$ has monomials of the form $b_{1}x^{\beta}y\frac{\partial}{\partial y}$ and/or $b_{\alpha + 1}y^{\alpha +1}\frac{\partial}{\partial y}$. Finally, in Case III we set $a_{\alpha + 1} = b_{\alpha + 1} = b_{0} = 0$, and it holds that $a_{\alpha} \neq 0$ and $a_{0}^{2} + b_{1}^{2} \neq 0$. Here, the component $P_{0}^{U}$ has monomials of the form $a_{\alpha}y^{\alpha}\frac{\partial}{\partial x}$ and/or $a_{0}x^{\beta }\frac{\partial}{\partial x}$, and $Q_{0}^{U}$ has a monomial of the form $b_{1}x^{\beta-1}y\frac{\partial}{\partial y}$.

All of these properties hold due to the vertices of $\mathcal{P}^{U}$ in each case. We adopt the notation in \eqref{eq-expression-pq} in order to have a unified approach and avoid the study of several analogous cases. Finally, the numbers $N_{1,2}$ in \eqref{eq-expression-pq} are positive integers that represent the highest ordinate of the vertices of $\Gamma_{1}^{(1)}$. In this proof we will work out the coefficients $A_{k} = A_{k}(\lambda)$. The other coefficients can be obtained with completely analogous computations.

Using \eqref{eq-expression-pq}, one obtains
\begin{align*}
&  Q_{0}^{U} (1,z  + \lambda)  - \displaystyle\frac{\beta}{\alpha}(z + \lambda)P_{0}^{U}(1,z + \lambda)  = \\ 
  & = \displaystyle\sum_{n = 0}^{\alpha + 1}b_{n}(z+\lambda)^{n} - \frac{\beta}{\alpha}\left(z\sum_{n = 0}^{\alpha + 1}a_{n}(z+\lambda)^{n} + \lambda\sum_{n = 0}^{\alpha + 1}a_{n}(z+\lambda)^{n} \right) \\
     & = \displaystyle\sum_{n = 0}^{\alpha + 1}\sum_{j = 0}^{n}\binom{n}{j}b_{n}\lambda^{n-j}z^{j} - \displaystyle\frac{\beta}{\alpha}\left(z\displaystyle\sum_{n = 0}^{\alpha + 1}\sum_{j = 0}^{n}\binom{n}{j}a_{n}\lambda^{n-j}z^{j} + \lambda \displaystyle\sum_{n = 0}^{\alpha + 1}\sum_{j = 0}^{n}\binom{n}{j}a_{n}\lambda^{n-j}z^{j}\right) \\
     & = \left( Q_{0}^{U}(1,\lambda) - \displaystyle\frac{\beta}{\alpha}\lambda P_{0}^{U}(1,\lambda) \right) \\
     & + \sum_{k = 1}^{\alpha + 1} \frac{z^{k}}{k!}\left( \sum_{n = k}^{\alpha + 1}\frac{n!}{(n-k)!}b_{n}\lambda^{n-k} - \frac{\beta}{\alpha}\left( \sum_{n = k-1}^{\alpha + 1}\frac{n!}{(n-k+1)!}a_{n}\lambda^{n-k+1} + \lambda\sum_{n = k}^{\alpha + 1}\frac{n!}{(n-k)!}a_{n}\lambda^{n-k}  \right) \right).
\end{align*}

The conclusion follows noticing that
$$\displaystyle\frac{\partial^{k}P_{0}^{U}}{\partial y^{k}}(1,\lambda) = \sum_{n = k}^{\alpha + 1}\frac{n!}{(n-k)!}a_{n}\lambda^{n-k}, \quad \quad \displaystyle\frac{\partial^{k}Q_{0}^{U}}{\partial y^{k}}(1,\lambda) = \sum_{n = k}^{\alpha + 1}\frac{n!}{(n-k)!}b_{n}\lambda^{n-k}.$$

Observe that we used the fact that $ Q_{0}^{U}(1,\lambda) - \displaystyle\frac{\beta}{\alpha}\lambda P_{0}^{U}(1,\lambda) = 0$.


\begin{thebibliography}{50}

\bibitem{AGR2011} A. Algaba, C. García, M. Reyes. \textit{Characterization of a monodromic singular point of a planar vector field}. \textbf{Nonlinear Analysis: Theory, Methods and Applications} 74(16) (2011), 5402--5414.

\bibitem{AGR2014} A. Algaba, C. García, M. Reyes. \textit{A new algorithm for determining the monodromy of a planar differential system}. \textbf{Applied Mathematics and Computation} 237 (2014), 419--429.

\bibitem{Alo2015} C. Alonso–Gonzalez. \textit{Infinitesimal Hartman–Grobman theorem in dimension three}. \textbf{An. Acad. Bras. Ciênc.} 87(3) (2015) 1499--1503.

\bibitem{AP90} D.K. Arrowsmith and C.M. Place. \textit{An introduction to Dynamical Systems}. Cambridge University Press (1990).

\bibitem{Ber78} F.S. Berezovskaya. \textit{Topological normal form for a system of two differential equations}. \textbf{Russian Math. Surv.} 33(2) (1978), 227--228.

\bibitem{BivHua19} C. Bivià-Ausina, J.A.C. Huarcaya. \textit{Polynomial maps with maximal multiplicity and the special closure}. \textbf{Monatsh Math} 118 (2019), 413--429.

\bibitem{BM90} M. Brunella, M. Miari. \textit{Topological equivalence of a plane vector field with its principal part defined through Newton Polyhedra}. \textbf{J. Diff. Equations} 85 (1990), 338--366.

\bibitem{Bru89} A.D. Bruno. \textit{Local methods in nonlinear differential equations}. Springer-Verlag Berlin Heidelberg (1989).

\bibitem{Bru2000} A.D. Bruno. \textit{Power Geometry in Algebraic and Differential Equations}. North-Holland Mathematical Library, 57. Amsterdam (2000).

\bibitem{DOP24} T.M. Dalbelo, R. Oliveira, O.H. Perez. \textit{Topological equivalence at infinity of a planar vector field and its principal part defined through Newton polytope}. \textbf{J. Diff. Equations} 408 (2024), 230--253.

\bibitem{DGV22} M.V. Demina, J. Giné, C. Valls. \textit{Puiseux integrability of differential equations}. \textbf{Qual. Theory Dyn. Syst.} 21 (2022), 35.


\bibitem{Dum77} F. Dumortier. \textit{Singularities of vector fields on the plane}. \textbf{J. Diff. Equations} 23 (1977), 53--106.

\bibitem{DLA06} F. Dumortier, J. Llibre, J.C. Artés. \textit{Qualitative theory of planar differential systems}. Universitext, Springer-Verlag Berlin Heidelberg (2006).

\bibitem{GH83} J. Guckenheimer, P. Holmes. \textit{Nonlinear Oscillations, Dynamical Systems, and Bifurcations of Vector Fields}. Springer-Verlag, New York 
(1983).

\bibitem{Kho77} A.G. Khovanskii. \textit{Newton polyhedra and toroidal varieties}. \textbf{Funkcional. Anal. i Priložen.} 11(4) (1977), 56--64.

\bibitem{Kho78} A.G. Khovanskii. \textit{Newton polyhedra and the genus of complete intersections}. \textbf{Funkcional. Anal. i Priložen.} 12(1) (1978), 51--61.

\bibitem{Kou76} A.G. Kouchnirenko. \textit{Polyèdres de Newton et nombres de Milnor}. \textbf{Invent Math} 32 (1976), 1--31.

\bibitem{OliVal26} R. Oliveira, C. Valls. \textit{Topological equivalence at infinity of second order planar vector fields and its principal part via Newton polytope}. To appear in \textbf{Israel Journal of Mathematics}. Preprint available at SSRN: http://dx.doi.org/10.2139/ssrn.5160056

\bibitem{Pan06} D. Panazzolo. \textit{Resolution of singularities of real-analytic vector fields in dimension three}. \textbf{Acta Math} 197(2) (2006), 167--289.

\bibitem{Pel95} M. Pelletier. \textit{\'Eclatements quasi homog\`enes}. \textbf{Ann. Fac. Sci. Toulouse Math} 4 (1995), 879--937.

\bibitem{PerSil22-ii} O.H. Perez, P.R. Silva. \textit{Resolution of singularities of 2-dimensional real analytic constrained differential systems}. \textbf{Bull. Sci. Math.} 179 (2022) 103179.

\bibitem{PerSil22} O.H. Perez, P.R. Silva. \textit{Singular impasse points of planar constrained differential systems}. \textbf{Bull. Belg. Math. Soc. Simon Stevin} 29 (2022) 611--643.

\bibitem{Tak74} F. Takens. \textit{Singularities of vector fields}. \textbf{Publ. Math. I.H.E.S.} 43 (1974), 47--100.

\bibitem{Zup96} V. \v{Z}upanovi\'c. \textit{Generalized principal part of some planar vector fields}. \textbf{Glasnik Matematicki} 31(51) (1996), 333--351.

\bibitem{Zup00} V. \v{Z}upanovi\'c. \textit{Topological equivalence of planar vector fields and their generalized principal part}. \textbf{J. Diff. Equations} 167 (2000), 1--15.


\end{thebibliography}
\end{document}